\documentclass[12pt,reqno]{amsart}
\usepackage{amsthm,amsfonts,amssymb,euscript}

\newcommand{\R}{\ensuremath{\mathbb{R}}}
\newcommand{\Z}{\ensuremath{\mathbb{Z}}}

\newcommand{\eps}{\ensuremath{\varepsilon}}

\newcommand{\E}{\ensuremath{\mathbf{E}}}
\newcommand{\F}{\ensuremath{\mathbf{F}}}
\newcommand{\N}{\ensuremath{\mathbf{N}}}
\newcommand{\B}{\ensuremath{\mathbf{B}}}

\newtheorem{theorem}{Theorem}[section]
\newtheorem{lemma}[theorem]{Lemma}
\newtheorem{proposition}[theorem]{Proposition}

\numberwithin{equation}{section}
\begin{document}

\title[Global well-posedness of the KP-I equation]{Global well-posedness of the KP-I initial-value problem in the energy space}
\author{A. D. Ionescu}
\address{University of Wisconsin--Madison}
\email{ionescu@math.wisc.edu}
\author{C. E. Kenig}
\address{University of Chicago}
\email{cek@math.uchicago.edu}
\author{D. Tataru}
\address{University of California--Berkeley}
\email{tataru@math.berkeley.edu}
\thanks{The first author was
supported in part by an NSF grant and a Packard Fellowship. The
second author was supported in part by an NSF grant. The
third author was supported in part by an NSF grant.}
\begin{abstract}
We prove that the KP-I initial-value problem
\begin{equation*}
\begin{cases}
\partial_tu+\partial_x^3u-\partial_x^{-1}\partial_y^2u+\partial_x(u^2/2)=0\text{ on }\mathbb{R}^2_{x,y}\times\mathbb{R}_t;\\
u(0)=\phi,
\end{cases}
\end{equation*}
is globally well-posed in the energy space
\begin{equation*}
\E^1(\R^2)=\{\phi:\mathbb{R}^2\to\mathbb{R}:\,\|\phi\|_{\mathbf{E}^1(\R^2)}\approx \|\phi\|_{L^2}+\|\partial_x\phi\|_{L^2}+\|\partial_x^{-1}\partial_y\phi\|_{L^2}<\infty\}.
\end{equation*}
\end{abstract}
\maketitle \tableofcontents

\section{Introduction}\label{section1}

In this paper we consider the KP-I initial-value problem
\begin{equation}\label{eq-1}
\begin{cases}
\partial_tu+\partial_x^3u-\partial_x^{-1}\partial_y^2u+\partial_x(u^2/2)=0;\\
u(0)=\phi,
\end{cases}
\end{equation}
on $\mathbb{R}^2_{x,y}\times\mathbb{R}_t$. The KP-I equation and
the KP-II equation, in which the sign of the term
$\partial_x^{-1}\partial_y^2u$ in \eqref{eq-1} is $+$ instead of
$-$, arise in physical contexts as models for the propagation of
dispersive long waves with weak transverse effects.

The KP-II equation is well understood from the point of view of
well-posedness: the KP-II initial-value problem is globally well-posed
for suitable data in $L^2$, on both $\mathbb{R}^2$ and
$\mathbb{T}^2=\mathbb{S}^1\times\mathbb{S}^1$, see \cite{Bo}, as well
as in some spaces larger than $L^2$, see \cite{TaTz} and the
references therein.

On the other hand, it has been shown in
\cite{MoSaTz1} that the KP-I initial-value problem is badly behaved
with respect to Picard iterative methods in standard Sobolev spaces,
since the flow map fails to be real-analytic at the origin in these
spaces.\footnote{Picard iterative methods can be applied, however, to
  produce local in time solutions for small low-regularity data in
  suitably weighted spaces, see \cite{CoIoKeSt}.} On the positive
side, it is known that the KP-I initial value problem is globally
well-posed in the ``second'' energy spaces on both $\mathbb{R}^2$ (see
\cite{Ke}, and also \cite{MoSaTz2} and \cite{MoSaTz3}) and
$\mathbb{T}^2$ (see \cite{IoKe}), as well as locally well-posed in
larger spaces.  These global well-posedness results rely on refined
energy methods. In this paper we show that the KP-I initial-value problem is globally
well-posed in the natural energy space of the equation.

Let $\xi$, $\mu$ and $\tau$ denote the Fourier variables with
respect to $x$, $y$ and $t$ respectively.  For $\sigma=1,2,\ldots$ we define the
Banach spaces $\mathbf{E}^\sigma=\mathbf{E}^\sigma(\R^2)$,
\begin{equation}\label{eq-2}
\mathbf{E}^\sigma
=\{\phi:\mathbb{R}^2\to\mathbb{R}:\,\|\phi\|_{\mathbf{E}^\sigma
}=\|\widehat{\phi}(\xi,\mu)\cdot
p(\xi,\mu)(1+|\xi|)^{\sigma}\|_{L^2_{\xi,\mu}}<\infty\},
\end{equation}
where $\widehat{\phi}$ denotes the Fourier transform of $\phi$ and
\begin{equation}\label{eq-3}
p(\xi,\mu)=1+\frac{|\mu|}{|\xi|+|\xi|^2}.
\end{equation}
Clearly,
\[
p(\xi,\mu)(1+|\xi|)^{\sigma}=(1+|\xi| )^\sigma +|\mu/ \xi|\cdot
(1+|\xi| )^{\sigma-1}.
\]
 Let
\[
\E^\infty=\bigcap_{\sigma=1}^\infty \E^\sigma
\]
with the induced metric. We recall the KP-I conservation laws (see,
for example, \cite{MoSaTz2} for formal justifications): if $t_1<t_2\in\mathbb{R}$
$u\in C([t_1,t_2]:\E^\infty)$ is a solution of the equation $\partial_tu+\partial_x^3u-\partial_x^{-1}\partial_y^2u+\partial_x(u^2/2)=0$
on $\mathbb{R}^2\times(t_1,t_2)$ then
\begin{equation}\label{conserve}
\widetilde{E}^0(u(t_1))=\widetilde{E}^0(u(t_2))\text{ and
}\widetilde{E}^1(u(t_1))=\widetilde{E}^1(u(t_2)),
\end{equation}
where, for any $\phi\in \E^1$,
\begin{equation}\label{E30}
\widetilde{E}^0(\phi)=\int_{\R^2} \phi^2\,dxdy,
\end{equation}
and
\begin{equation}\label{E31}
\widetilde{E}^1(\phi)=\int_{\R^2}(\partial_x\phi
)^2\,dxdy+\int_{\R^2}(\partial_x^{-1}\partial_y\phi
)^2\,dxdy-\frac{1}{3}\int_{\R^2} \phi^3\,dxdy.
\end{equation}
Consequently, if $t_0\in[t_1,t_2]$ and $\|u(t_0)\|_{\E^1}\leq 1$ then we have the uniform bound
\begin{equation}\label{conserve2}
\sup_{t\in[t_1,t_2]}\|u(t)\|_{\E^1}\lesssim \|u(t_0)\|_{\E^1}.
\end{equation}
Our main theorem concerns global well-posedness of the KP-I
initial-value problem in the energy space $\E^1$.

\newtheorem{Main1}{Theorem}[section]
\begin{Main1}\label{Main1}
(a) Assume $\phi\in \E^\infty$. Then there is a unique global
solution
\begin{equation*}
u=S^\infty(\phi)\in C(\R:\E^\infty)
\end{equation*}
of the initial-value problem \eqref{eq-1}. In addition, for any
$T\in[0,\infty)$ and any $\sigma\in\{1,2,\ldots\}$
\begin{equation}\label{ww1}
\sup_{|t|\leq T}\|S^\infty(\phi)(t)\|_{\E^\sigma}\leq C(T,\sigma,\|\phi\|_{\E^\sigma}).
\end{equation}

(b) Assume $T\in\R_+$. Then the mapping
\[
S^\infty_T=\mathbf{1}_{[-T,T]}(t)\cdot S^\infty:\E^\infty\to
C([-T,T]:\E^\infty)
\]
 extends uniquely to a continuous mapping
\[
S^1_T:\E^1\to C([-T,T]:\E^1),
\]
 and
\begin{equation*}
\widetilde{E}^j(u(t))=\widetilde{E}^j(\phi)\text{ for any
}t\in[-T,T] \text{ and }j\in\{0,1\}.
\end{equation*}
\end{Main1}

We remark that the global existence of smooth solutions stated in Theorem \ref{Main1} (a) is also new. The earlier global existence theorems of \cite{Ke}, \cite{MoSaTz2}, and \cite{MoSaTz3} rely of the conservation of the second energy, which requires the stronger momentum condition $\partial_x^{-2}\partial_y^2\phi\in L^2$. Also, a simple additional argument shows that \eqref{ww1} can be improved to
\begin{equation*}
\sup_{|t|\leq T}\|S^\infty(\phi)(t)\|_{\E^\sigma}\leq C(T,\sigma,\|\phi\|_{\E^1})\cdot \|\phi\|_{\E^\sigma}.
\end{equation*}

We discuss now some of the ingredients in the proof of Theorem
\ref{Main1}.  One might try a direct perturbative approach (which goes
back to work on the KdV equation in \cite{KePoVe2}, \cite{Bo2}, \cite{KePoVe}, and
nonlinear wave equations in \cite{KlMa}), based on the properties of
solutions to the linear equation
\begin{equation}
\label{eq-lin}
\begin{cases}
\partial_tu+\partial_x^3u-\partial_x^{-1}\partial_y^2u=f;\\
u(0)=\phi.
\end{cases}
\end{equation}
For some suitable spaces $\F^1(T)$ and $\N^1(T)$
one would like to prove a linear bound for solutions
to \eqref{eq-lin} on $\R^2\times[-T,T]$,
$T\in(0,1]$, of the form
\begin{equation}\label{ch1a}
\|u\|_{\F^1(T)}\lesssim \|\phi\|_{\E^1}+\|f\|_{\N^1(T)},\\
\end{equation}
together with a matching nonlinear estimate
\begin{equation}\label{ch1b}
\|-\partial_x(u^2/2)\|_{\N^1(T)}\lesssim  \|u\|_{\F^1(T)}^2.
\end{equation}
Due to \cite{MoSaTz1}, it is known however that the inequalities
\eqref{ch1a}, \eqref{ch1b} cannot hold for any choice of the spaces
$\F^1(T)$ and $\N^1(T)$; this forces us to approach the problem in a
less perturbative way.

To prove Theorem \ref{Main1} (a) we define instead the normed spaces
$\F^1(T)$, $\N^1(T)$, and the semi-normed space $\B^1(T)$ and show
that if $u$ is a smooth solution of \eqref{eq-1} on
$\R^2\times[-T,T]$, $T\in(0,1]$, then
\begin{equation}\label{ch2}
\begin{cases}
&\|u\|_{\F^1(T)}\lesssim \|u\|_{\B^1(T)}+\|-\partial_x(u^2/2)\|_{\N^1(T)};\\
&\|-\partial_x(u^2/2)\|_{\N^1(T)}\lesssim \|u\|_{\F^1(T)}^2;\\
&\|u\|_{\B^1(T)}^2\lesssim \|\phi\|_{\E^1}^2+\|u\|_{\F^1(T)}^3.
\end{cases}
\end{equation}
The inequalities \eqref{ch2} and a simple continuity argument still
suffice to control $\|u\|_{\F^1(T)}$, provided that
$\|\phi\|_{\E^1}\ll 1$ (which can be arranged by rescaling). The first
inequality in \eqref{ch2} is the analogue of the linear estimate
\eqref{ch1a}, and uses the linear equation \eqref{eq-lin}. The second
inequality in \eqref{ch2} is the analogue of the bilinear estimate
\eqref{ch1b}. The last inequality in \eqref{ch2} is an energy-type
estimate.

To prove Theorem \ref{Main1} (b) we need to  exploit
several special symmetries of the equation satisfied by the
difference of two solutions. This difference equation has special
symmetries for real-valued solutions in $L^2$ and in
$\dot{H}^{-1}_x$. To exploit these symmetries, we define the
normed spaces $\overline{\F}^0$, $\overline{\N}^0$, and the
semi-normed space $\overline{\B}^0$, and prove a second set of
linear, bilinear, and energy estimates, similar to \eqref{ch2}.
Then we adapt the Bona-Smith method \cite{BoSm} to prove the
continuity of the flow in the space $\E^1$.

We explain now our strategy to define the main normed and
semi-normed spaces. Ideally, one would like to use standard
$X^{s,b}$- type structures (as in \cite{Bo2}, \cite{KePoVe}) for
the spaces $\F^1(T)$ and $\N^1(T)$. For such spaces, however, the
bilinear estimate $\|\partial_x(uv)\|_{\N^1(T)}\lesssim \|u\|_{\F^1(T)}\|v\|_{\F^1(T)}$ cannot hold even for solutions
$u,v$ of the linear homogeneous equation. This bilinear estimate
is only possible if we weaken significantly the contributions of
the components of the functions $u$ and $v$ of high frequency and
low modulation. To achieve this we still use $X^{s,b}$-type
structures for the spaces $\F^1(T)$ and $\N^1(T)$, but only on
small, frequency dependent time intervals. A similar method was
used recently in \cite{ChCoTa} and \cite{KoTa} to prove a-priori
bounds for the $1$-d cubic nonlinear Schr\"{o}dinger equation in
negative Sobolev spaces.

The second step is to define
$\|u\|_{\B^1(T)}$ sufficiently large to be able to still prove the
linear estimate $\|u\|_{\F^1(T)}\lesssim \|u\|_{\B^1(T)}+\|-\partial_x(u^2/2)\|_{\N^1(T)}$. Finally, we
use frequency-localized energy estimates and the symmetries of the
equation \eqref{eq-1}\footnote{The two main symmetries used at
this stage are the fact that the solution $u$ is real-valued and
the precise form of the nonlinearity $-\partial_x(u^2/2)$.} to
prove the energy estimate $\|u\|_{\B^1(T)}^2\lesssim\|\phi\|_{\E^1}^2+\|u\|_{\F^1(T)}^3$. These symmetries allow us to trade high frequencies
for low frequencies in trilinear forms, improving the timescale from
frequency dependent time intervals (as guaranteed by the bilinear
estimates) to frequency independent time intervals.

A new twist arises in the proof of part (b) of Theorem \ref{Main1}. The symmetries of the difference equation are not as good as the symmetries of the nonlinear equation, which causes difficulties in the proofs of suitable energy estimates.
The low frequency part of the solution turns out to be particularly
harmful in the difference equation. To avoid this difficulty  we
define the normed spaces $\overline{\F}^0$, $\overline{\N}^0$, and the
semi-normed space $\overline{\B}^0$, which have a special low-frequency
structure.

The rest of the paper is organized as follows: in section
\ref{definitions} we summarize most of the notation, define the
main normed spaces, and prove some of their basic properties. In
section \ref{glolibien} we state our main global linear,
bilinear, and energy estimates. The proof of the bilinear estimate
Proposition \ref{Lemmab1} depends on the dyadic bilinear estimates
proved in sections \ref{dyadic1} and
\ref{dyadic3}; the energy estimates Proposition \ref{Lemmad1} and
Proposition \ref{Lemmak4} are proved in section \ref{energyproof}.
In section \ref{proof1} we prove the main theorem, using the
linear, bilinear, and energy estimates of section \ref{glolibien}.
In section \ref{L2bi}, which is self-contained, we prove the
bilinear $L^2$ estimates in Corollary \ref{Main9co}; these $L^2$
estimates are the main building blocks in all the dyadic estimates
in sections \ref{energyproof}, \ref{dyadic1}, and
\ref{dyadic3}. The key technical ingredient is the scale-invariant
estimate in Lemma \ref{Main9} (a), which is also used in
\cite{CoIoKeSt}. In section \ref{energyproof} we prove the energy
estimates Proposition \ref{Lemmad1} and Proposition \ref{Lemmak4}.
Finally, in sections \ref{dyadic1} and
\ref{dyadic3} we prove the dyadic bilinear estimates used in
Proposition \ref{Lemmab1}.

\section{Notation and definitions}\label{definitions}

Let $\mathbb{Z}_+=\mathbb{Z}\cap[0,\infty)$. For $k\in \mathbb{Z}$ let
\begin{equation*}
I_k=\{\xi:|\xi|\in[(3/4)\cdot 2^k,(3/2)\cdot 2^k]\}\text{ and }\widetilde{I}_k=\{\xi:|\xi|\in[2^{k-1},2^{k+1}]\}.
\end{equation*}
Let $\eta_0:\mathbb{R}\to[0,1]$ denote an even smooth function
supported in $[-8/5,8/5]$ and equal to $1$ in $[-5/4,5/4]$. For
$k\in\mathbb{Z}\cap [1,\infty)$ let
$\eta_k(\xi)=\eta_0(\xi/2^k)-\eta_0(\xi/2^{k-1})$. For $k\in\Z_+$
let $\eta_{\leq k}=\sum_{k'=0}^k\eta_{k'}$. For $k\in\mathbb{Z}$
let $\chi_k(\xi)=\eta_0(\xi/2^k)-\eta_0(\xi/2^{k-1})$. For
$(\xi,\mu)\in\mathbb{R}\setminus\{0\}\times\mathbb{R}$ let
\begin{equation}\label{omega}
\omega(\xi,\mu)=\xi^3+\mu^2/ \xi.
\end{equation}

For $k\in\Z$ we define the dyadic $X^{s,b}$-type normed spaces $X_k=X_k(\R^3)$,
\begin{equation}\label{sp7}
\begin{split}
X_{k}=\{f\in L^2&(\R^3):\,f\text{ is supported in }\widetilde{I}_k\times\R\times\R\text{ and }\\
&\|f\|_{X_k}=\sum_{j=0}^{\infty}2^{j/2}\|\eta_j(\tau-\omega(\xi,\mu)) \cdot f\|_{L^2}<\infty\}.
\end{split}
\end{equation}
We use an $l^1$ Besov-type norm with respect to modulations.
Structures of this type were introduced, for instance, in \cite{Ta},
and are useful in order to prevent high frequency losses in bilinear
and trilinear estimates.

The definition shows easily that if $k\in\Z$ and $f_k\in X_k$ then
\begin{equation}\label{sp7.1}
\left\|\int_\R
|f_k(\xi,\mu,\tau')|\,d\tau'\right\|_{L^2_{\xi,\mu}}\lesssim \|f_k\|_{X_k}.
\end{equation}
Moreover, if $k\in\Z$, $l\in\Z_+$, and $f_k\in X_k$ then
\begin{equation}\label{sp7.2}
\begin{split}
&\sum_{j=l+1}^\infty 2^{j/2}\left\|\eta_{j}(\tau-\omega(\xi,\mu))\cdot \int_\R |f_k(\xi,\mu,\tau')|\cdot 2^{-l}(1+2^{-l}|\tau-\tau'| )^{-4}\,d\tau'\right\|_{L^2}\\
&+2^{l/2}\left\|\eta_{\leq l}(\tau-\omega(\xi,\mu))\cdot \int_\R |f_k(\xi,\mu,\tau')|\cdot 2^{-l}(1+2^{-l}|\tau-\tau'| )^{-4}\,d\tau'\right\|_{L^2}\\
&\lesssim \|f_k\|_{X_k}.
\end{split}
\end{equation}
In particular, if $k\in\Z$, $l\in\Z_+$, $t_0\in\R$, $f_k\in X_k$, and $\gamma\in\mathcal{S}(\R)$, then
\begin{equation}\label{sp7.3}
\|\mathcal{F}[\gamma(2^l(t-t_0))\cdot \mathcal{F}^{-1}(f_k)]\|_{X_k}\lesssim \|f_k\|_{X_k}.
\end{equation}

For $k\in\Z$ let $k_+=\max(k,0)$, and let $P_k$ denote the
operator on $L^2(\mathbb{R}^3)$ defined by the Fourier multiplier
$(\xi,\mu,\tau)\to\mathbf{1}_{I_k}(\xi)$. By a slight abuse of
notation, we also let $P_k$ denote the operator on
$L^2(\mathbb{R}^2)$ defined by the Fourier multiplier
$(\xi,\mu)\to\mathbf{1}_{I_k}(\xi)$. For $l\in\Z$ let
\[
P_{\leq l}=\sum_{k\leq l}P_k, \qquad P_{\geq l}=\sum_{k\geq l}P_k.
\]
With $p$ as in \eqref{eq-3}, for $k\in\Z$ define the
frequency localized initial data spaces
\begin{equation}\label{gp1}
E_k=\{\phi:\R^2\to\R:\mathcal{F}(\phi)=\mathbf{1}_{\widetilde{I}_k}(\xi)\mathcal{F}(\phi)\text{ and
}\|\phi\|_{E_k}=\|\widehat{\phi}\cdot
p(\xi,\mu)\|_{L^2_{\xi,\mu}}<\infty\},
\end{equation}
and
\begin{equation}\label{gp1.1}
\overline{E}_k=\{\phi:\R^2\to\R:\mathcal{F}(\phi)=\mathbf{1}_{\widetilde{I}_k}(\xi)\mathcal{F}(\phi)\text{ and
}\|\phi\|_{\overline{E}_k}=\|\widehat{\phi}\|_{L^2_{\xi,\mu}}<\infty\}.
\end{equation}

The corresponding frequency localized energy spaces for the solutions are
\[
C_0(\R:E_k)=\{u_k\in C(\R:E_k):u_k\text{ is supported in
}\R^2\times[-4,4]\}
\]
\[
C_0(\R:\overline{E}_k)=\{u_k\in
C(\R:\overline{E}_k):u_k\text{ is supported in
}\R^2\times[-4,4]\}.
\]

At frequency $2^k$ we will use the $X^{s,b}$ structure
given by the $X_k$ norm, uniformly on the $2^{-k_+}$ time scale.
 For $k\in\Z$ we define the normed spaces
\begin{equation}\label{sp5.1}
\begin{split}
F_k=&\{u_k\in
C_0(\R:E_k):\,\|u_k\|_{F_k}=\sup_{t_k\in\R}\|p(\xi,\mu)\cdot
\mathcal{F}[u_k\cdot \eta_0(2^{k_+}(t-t_k))]\|_{X_k}<\infty\},
\end{split}
\end{equation}
and
\begin{equation}\label{sp5.1.1}
\overline{F}_k=\{u_k\in
C_0(\R:\overline{E}_k):\,\|u_k\|_{\overline{F}_k}=\sup_{t_k\in\R}\|\mathcal{F}[u_k\cdot
\eta_0(2^{k_+}(t-t_k))]\|_{X_k}<\infty\}.
\end{equation}
For $k\in\Z$ we define the normed spaces $N_k=C_0(\R:E_k)$ and
$\overline{N}_k=C_0(\R:\overline{E}_k)$ (as vector spaces), which are used to measure the
frequency $2^k$ part of the nonlinear term, with norms
\begin{equation}\label{gp2}
\|f_k\|_{N_k}=\sup_{t_k\in\R}\|p(\xi,\mu)(\tau-\omega(\xi,\mu)+i2^{k_+})^{-1}\cdot
\mathcal{F}[f_k\cdot \eta_0(2^{k_+}(t-t_k))]\|_{X_k},
\end{equation}
and
\begin{equation}\label{gp2.1}
\|u_k\|_{\overline{N}_k}=\sup_{t_k\in\R}\|(\tau-\omega(\xi,\mu)+i2^{k_+})^{-1}\cdot
\mathcal{F}[f_k\cdot \eta_0(2^{k_+}(t-t_k))]\|_{X_k}.
\end{equation}

The bounds we obtain for solutions of the KP-I equation are on a fixed
time interval, while the above function spaces are not.  To remedy
this, for any time $T\in(0,1]$ we define the normed spaces
\begin{equation}\label{sb30}
\begin{cases}
&F_k(T)=\{u_k\in
C([-T,T]:E_k):\|u_k\|_{F_k(T)}=\inf\limits_{\widetilde{u}_k=u_k\text{
in }\R^2
\times [-T,T]}\|\widetilde{u}_k\|_{F_k}<\infty\};\\
&N_k(T)=\{f_k\in
C([-T,T]:E_k):\|u_k\|_{N_k(T)}=\inf\limits_{\widetilde{f}_k=f_k\text{
in }\R^2\times [-T,T]}\|\widetilde{u}_k\|_{N_k}<\infty\},
\end{cases}
\end{equation}
where the infimum is taken over all extensions $\widetilde{u}_k\in
C_0(\R:E_k)$ of $u_k$. Similarly we define the normed spaces
\begin{equation}\label{sb30.1}
\begin{cases}
&\negmedspace\negmedspace\overline{F}_k(T)=\{u_k\in
C([-T,T]:\overline{E}_k):\|u_k\|_{\overline{F}_k(T)}=\inf\limits_{\widetilde{u}_k=u_k\text{
in }
\R^2\times [-T,T]}\|\widetilde{u}_k\|_{\overline{F}_k}<\infty\};\\
&\negmedspace\negmedspace\overline{N}_k(T)=\{f_k\in
C([-T,T]:\overline{E}_k):\|f_k\|_{\overline{N}_k(T)}=
\inf\limits_{\widetilde{f}_k=f_k\text{ in }\R^2\times
[-T,T]}\|\widetilde{f}_k\|_{\overline{N}_k}<\infty\},
\end{cases}
\end{equation}
where the infimum is taken over all extensions $\widetilde{u}_k\in
C_0(\R:\overline{E}_k)$ of $u_k$.

So far we have defined the dyadic function spaces where we measure the
solution $u$ to \eqref{eq-1} and the nonlinearity.  We assemble these
in a straightforward manner using a Littlewood-Paley decomposition to
obtain the global function spaces for the solutions.  In what follows
we let $\sigma\in\Z_+$ and $T \in (0,1]$. We define the Banach spaces
for the initial data
\begin{equation}\label{ty1}
\E^\sigma=\{\phi:\R^2\to\R:
\|\phi\|_{\mathbf{E}^\sigma}=\|\widehat{\phi}\cdot
p(\xi,\mu)(1+|\xi| )^\sigma\|_{L^2_{\xi,\mu}}<\infty\},
\end{equation}
and
\begin{equation}\label{ty1.1}
\overline{\E}^{\sigma} =\{\phi:\R^2\to\R:
\|\phi\|_{\overline{\mathbf{E}}^\sigma }=\|\widehat{\phi}\cdot
(1+|\xi|^{-5/8})(1+|\xi| )^\sigma\|_{L^2_{\xi,\mu}}<\infty\}.
\end{equation}
Their intersections are denoted by
\[
\E^\infty=\bigcap_{\sigma=1}^\infty\E^\sigma, \qquad
\overline{\E}^{\infty}=\bigcap_{\sigma=1}^\infty\overline{\E}^{\sigma}.
\]
For $u\in C([-T,T]:\E^\infty)$, respectively $u\in C([-T,T]:\overline{\E}^\infty)$, we define
\begin{equation}\label{sp6}
\|u\|_{\B^\sigma(T)}^2=\|P_{\leq 0}(u(0))\|^2_{\E^\sigma}
+\sum_{k\geq 1}\sup_{t_k\in[-T,T]}2^{2\sigma
k}\|P_k(u(t_k))\,\|_{E_k}^2,
\end{equation}
and
\begin{equation}\label{sp6.6.1}
\|u\|_{\overline{\B}^\sigma(T)}^2=\|P_{\leq
0}(u(0))\|^2_{\overline{\E}^{0}} +\sum_{k\geq
1}\sup_{t_k\in[-T,T]}2^{2\sigma
k} \|P_k(u(t_k))\,\|_{\overline{E}_k}^2.
\end{equation}
Notice that the $\overline{\E}^{\sigma}$ and
$\overline{\B}^\sigma(T)$ norms, which are used for the difference
equation, have the added factor $(1+|\xi|^{-5/8})$. This gives extra
decay at low frequencies, and is essential in our analysis.

Finally, the $X^{s,b}$- type control of the solutions, respectively
the nonlinearity is achieved using  the normed spaces
\begin{equation}\label{sp5}
\begin{split}
&\mathbf{F}^\sigma(T)= \{u\in C([-T,T]:\E^\infty):
\|u\|_{\mathbf{F}^\sigma(T)}^2=\sum_{k\in \Z}2^{2\sigma
k_+}\|P_k(u)\|_{F_k(T)}^2<\infty\},\\
&\mathbf{N}^\sigma(T)= \{u\in
C([-T,T]:\E^\infty):\|u\|_{\N^\sigma(T)}^2=\sum_{k\in\Z}2^{2\sigma
k_+}\|P_k(u)\|^2_{N_k(T)}<\infty\}.\\
\end{split}
\end{equation}
For the difference equation we use the normed spaces
\begin{equation}\label{sp5.11}
\begin{split}
&\overline{\mathbf{F}}^0(T)= \{u\in
C([-T,T]:\overline{\E}^{\infty}):
\negmedspace\|u\|_{\overline{\mathbf{F}}^0
(T)}^2\negmedspace=\negmedspace\sum_{k\in
\Z}(1+2^{-5k/4})\|P_k(u)\|_{\overline{F}_k(T)}^2<\infty\},\\
&\overline{\mathbf{N}}^0(T)= \{u\in
C([-T,T]:\overline{\E}^{\infty}):\negmedspace
\|u\|_{\overline{\N}^0(T)}^2\negmedspace=\negmedspace\sum_{k\in\Z}(1+2^{-5k/4})\
\|P_k(u)\|^2_{\overline{N}_k(T)}<\infty\}.
\end{split}
\end{equation}

For any $k\in\Z$ we define the set $S_k$ of {\it{k-acceptable time
multiplication factors}}
\begin{equation}\label{sb20}
S_k=\{m_k:\R\to\R:\|m_k\|_{S_k}=\sum_{j=0}^{10}2^{-jk_+}\|\partial^jm_k\|_{L^\infty}<\infty\}.
\end{equation}
Direct estimates using the definitions and  \eqref{sp7.2} show
that for any $\sigma\in\Z_+$ and $T\in(0,1]$
\begin{equation}\label{sb21}
\begin{cases}
&\left\|\sum_{k\in\Z}m_k(t)\cdot P_k(u)\right\|_{\F^\sigma(T)}\lesssim (\sup_{k\in\Z}\|m_k\|_{S_k})\cdot \|u\|_{\F^\sigma(T)};\\
&\left\|\sum_{k\in\Z}m_k(t)\cdot P_k(u)\right\|_{\N^\sigma(T)}\lesssim (\sup_{k\in\Z}\|m_k\|_{S_k})\cdot \|u\|_{\N^\sigma(T)};\\
&\left\|\sum_{k\in\Z}m_k(t)\cdot
P_k(u)\right\|_{\B^\sigma(T)}\lesssim
(\sup_{k\in\Z}\|m_k\|_{S_k})\cdot \|u\|_{\B^\sigma(T)},
\end{cases}
\end{equation}
and
\begin{equation}\label{sb21.1}
\begin{cases}
&\left\|\sum_{k\in\Z}m_k(t)\cdot
P_k(u)\right\|_{\overline{\F}^0 (T)}
\lesssim (\sup_{k\in\Z}\|m_k\|_{S_k})\cdot \|u\|_{\overline{\F}^0 (T)};\\
&\left\|\sum_{k\in\Z}m_k(t)\cdot P_k(u)\right\|_{\overline{\N}^0 (T)}
\lesssim (\sup_{k\in\Z}\|m_k\|_{S_k})\cdot \|u\|_{\overline{\N}^0 (T)};\\
&\left\|\sum_{k\in\Z}m_k(t)\cdot
P_k(u)\right\|_{\overline{\B}^0 (T)}\lesssim
(\sup_{k\in\Z}\|m_k\|_{S_k})\cdot \|u\|_{\overline{\B}^0 (T)}.
\end{cases}
\end{equation}

\section{Global linear, bilinear and energy
estimates}\label{glolibien}

In this section we state our main linear, bilinear and energy
estimates. We show first that $\F^\sigma(T)\hookrightarrow
C([-T,T]:\E^\sigma)$ for $\sigma\in\Z_+$, $T\in(0,1]$.

\begin{lemma}\label{Lemmaa2}
If $\sigma\in\Z_+$, $T\in(0,1]$, and $u\in \F^\sigma(T)$, then
\begin{equation}\label{hh80}
\sup_{t\in[-T,T]}\|u(t)\|_{\E^\sigma}\lesssim  \|u\|_{\F^\sigma(T)}.
\end{equation}
\end{lemma}

\begin{proof}[Proof of Lemma \ref{Lemmaa2}] In view of the definitions,
it suffices to prove that if $k\in\Z$, $t_k\in[-1,1]$, and
$\widetilde{u}_k\in F_k$ then
\begin{equation}\label{sb1}
\|p(\xi,\mu)\cdot
\mathcal{F}[\widetilde{u}_k(t_k)]\|_{L^2_{\xi,\mu}}\lesssim
\|p(\xi,\mu)\cdot \mathcal{F}[\widetilde{u}_k\cdot
\eta_0(2^{k_+}(t-t_k))]\|_{X_k}.
\end{equation}
Let $f_k=\mathcal{F}[\widetilde{u}_k\cdot
\eta_0(2^{k_+}(t-t_k))]$, so
$$\mathcal{F}[\widetilde{u}_k(t_k)](\xi,\mu)=c\int_\R
f_k(\xi,\mu,\tau)e^{it_k\tau}\,d\tau.$$ Thus, using \eqref{sp7.1},
\begin{equation*}
\|p(\xi,\mu)\cdot
\mathcal{F}[\widetilde{u}_k(t_k)]\|_{L^2_{\xi,\mu}}\lesssim
\|p(\xi,\mu)\cdot f_k\|_{L^2_{\xi,\mu}L^1_\tau}\lesssim
\|p(\xi,\mu)\cdot f_k\|_{X_k},
\end{equation*}
which gives \eqref{sb1}.
\end{proof}

We prove now a linear estimate.

\begin{proposition}\label{Lemmaa3}
Assume $T\in(0,1]$, $u,v\in C([-T,T):\E^\infty)$ and
\begin{equation}\label{sb2a}
\partial_t u+\partial_x^3u-\partial_x^{-1}\partial_y^2u=v\text{ on }\R^2\times(-T,T).
\end{equation}

(a) Then, for any $\sigma\in\Z_+$,
\begin{equation}\label{amy20a}
\|u\|_{\F^\sigma(T)}\lesssim  \|
u\|_{\B^\sigma(T)}+\|v\|_{\N^\sigma(T)}.
\end{equation}

(b) Assume, in addition, that $u(0)\in \overline{\E}^\infty$ and
$v\in\overline{\N}^0(T)$. Then $u\in \overline{\F}^0(T)$ and
\begin{equation}\label{amy20b}
\|u\|_{\overline{\F}^0(T)}\lesssim
\|u\|_{\overline{\B}^0(T)}+\|v\|_{\overline{\N}^0(T)}.
\end{equation}
\end{proposition}

\begin{proof}[Proof of Proposition \ref{Lemmaa3}] In view of the definitions,
it suffices to prove that if $k\in\Z$ and $u,v\in C([-T,T]:E_k)$
solve \eqref{sb2a}, then
\begin{equation}\label{amy20}
\begin{cases}
&\|P_k(u)\|_{F_k(T)}\lesssim \|P_k(u(0))\|_{E_k}+\|v\|_{N_k(T)}\text{ if }k\leq 0;\\\\
&\|P_k(u)\|_{F_k(T)}\lesssim
\sup\limits_{t_k\in[-T,T]}\|P_k(u(t_k))\|_{E_k}+\|v\|_{N_k(T)} \text{
if }k\geq 1.
\end{cases}
\end{equation}
Let $\widetilde{v}\in C_0(\R:E_k)$ denote an extension of
  $v$ such that $\|\widetilde{v}\|_{N_k}\leq C\|v\|_{N_k(T)}$. Using
  \eqref{sb21}, we may assume that $\widetilde{v}$ is supported in
  $\mathbb{R}^2\times [-T-2^{-k_+-10},T+2^{-k_+-10}]$, $k\in\Z$. For
  $t\in\R$ let $W(t)$ denote the solution at time $t$ of the free KP-I
  evolution, i.e. the operator on $L^2(\R^2)$ defined by the Fourier
  multiplier $(\xi,\mu)\to e^{it\omega(\xi,\mu)}$. For $t\geq T$ we
  define
\begin{equation*}
\widetilde{u}(t)=\eta_0(2^{k_++5}(t-T))\Big[W(t-T)P_k(u(T))+\int_{T}^tW(t-s)(P_k(\widetilde{v}(s)))\,ds\Big].
\end{equation*}
For $t\leq -T$ we define
\begin{equation*}
\widetilde{u}(t)=\eta_0(2^{k_++5}(t+T))\Big[W(t+T)P_k(u(-T))+\int_{-T}^tW(t-s)(P_k(\widetilde{v}(s)))\,ds\Big].
\end{equation*}
With $\widetilde{u}(t)=u(t)$ for $t\in[-T,T]$, it is clear that
$\widetilde{u}\in C_0(\R:E_k)$ is an extension of $u$. Also, using
\eqref{sb21},
\begin{equation*}
\|u\|_{F_k(T)}\lesssim \sup_{t_k\in[-T,T]}\|p(\xi,\mu)\cdot
\mathcal{F}[\widetilde{u}\cdot \eta_0(2^{k_+}(t-t_k))]\|_{X_k},
\end{equation*}
where the supremum is taken over $t_k\in [-T,T]$.

In view of the definitions and \eqref{sb21}, it suffices to prove
that if $k\in\Z$, $\phi_k\in E_k$, and $v_k\in N_k$ then
\begin{equation}\label{sb4}
\begin{split}
\|p(\xi,\mu)\cdot \mathcal{F}[u_k\cdot
&\eta_0(2^{k_+}t)]\|_{X_k}\lesssim \|p(\xi,\mu)\cdot \widehat{\phi_k}\|_{L^2_{\xi,\mu}}\\
&+\|p(\xi,\mu)(\tau-\omega(\xi,\mu)+i2^{k_+})^{-1}\cdot
\mathcal{F}(v_k)\|_{X_k},
\end{split}
\end{equation}
where
\begin{equation}\label{sb5}
u_k(t)=W(t)(\phi_k)+\int_{0}^tW(t-s)(v_k(s))\,ds.
\end{equation}
It follows from \eqref{sb5} that
\begin{equation*}
\begin{split}
\mathcal{F}&[u_k\cdot \eta_0(2^{k_+}t)](\xi,\mu,\tau)=\widehat{\phi}_k(\xi,\mu)\cdot 2^{-k_+}\widehat{\eta_0}(2^{-k_+}(\tau-\omega(\xi,\mu)))\\
&+C\int_\R\mathcal{F}(v_k)(\xi,\mu,\tau')\cdot
\frac{2^{-k_+}\widehat{\eta_0}(2^{-k_+}(\tau-\tau'))-2^{-k_+}\widehat{\eta_0}(2^{-k_+}(\tau-\omega(\xi,\mu)))}{\tau'-\omega(\xi,\mu)}\,d\tau'.
\end{split}
\end{equation*}
We observe now that
\begin{equation*}
\begin{split}
\Big|&\frac{2^{-k_+}\widehat{\eta_0}(2^{-k_+}(\tau-\tau'))-2^{-k_+}\widehat{\eta_0}(2^{-k_+}(\tau-\omega(\xi,\mu)))}{\tau'-\omega(\xi,\mu)}\cdot (\tau-\omega(\xi,\mu)+i2^{k_+})\Big|\\
&\lesssim 2^{-k_+}(1+2^{-k_+}|\tau-\tau'|
)^{-4}+2^{-k_+}(1+2^{-k_+}|\tau-\omega(\xi,\mu)| )^{-4},
\end{split}
\end{equation*}
and the bound \eqref{sb4} follows from \eqref{sp7.1} and
\eqref{sp7.2}.
\end{proof}

We continue with our main bilinear estimates.

\begin{proposition}\label{Lemmab1}
a) If $\sigma\in\{1,2,3\}$, $T\in[0,1)$, and $u,v\in \F^\sigma(T)$
then
\begin{equation}\label{on1}
\|\partial_x(uv)\|_{\N^\sigma(T)}\lesssim
\|u\|_{\F^\sigma(T)}\cdot \|v\|_{\F^1(T)}+\|u\|_{\F^1(T)}\cdot
\|v\|_{\F^\sigma(T)}.
\end{equation}

b) If $T\in(0,1]$,  $u \in \overline{\F}^0(T)$ and $v \in \F^1(T)$
then
\begin{equation}\label{on2}
\|\partial_x(uv)\|_{\overline{\N}^0(T)}\lesssim
\|u\|_{\overline{\F}^0(T)}\cdot \|v\|_{\F^1(T)}.
\end{equation}
\end{proposition}

\begin{proof}[Proof of Proposition \ref{Lemmab1}] We fix
extensions $\widetilde{u},\widetilde{v}\in C_0(\R:\E^\infty)$ of
$u,v$ such that $\|P_k(\widetilde{u})\|_{F_k}\leq
2\|P_k(u)\|_{F_k(T)}$ and $\|P_k(\widetilde{v})\|_{F_k}\leq
2\|P_k(v)\|_{F_k(T)}$, $k\in\Z$. Let
$\widetilde{u}_k=P_k(\widetilde{u})$ and
$\widetilde{v}_k=P_k(\widetilde{v})$, $k\in\Z$. It follows from
Lemma \ref{Lemmab2}, Lemma \ref{Lemmab3}, and Lemma \ref{Lemmab6}
that
\begin{equation}\label{th1}
\begin{cases}
&2^{k_+}\|P_k(\partial_x(\widetilde{u}_{k_1}\widetilde{v}_{k_2}))\|_{N_k}\lesssim
2^{-|k_1|/2}
2^{{k_1}_+}\|\widetilde{u}_{k_1}\|_{F_{k_1}}\cdot
2^{{k_2}_+}\|\widetilde{v}_{k_2}\|_{F_{k_2}}\\
&\text{if }k_1\leq k_2\text{ and }|k_2-k|\leq 40.
\end{cases}
\end{equation}
Also, it follows from Lemma \ref{Lemmab4}, Lemma \ref{Lemmab5},
and Lemma \ref{Lemmab6} that
\begin{equation}\label{th2}
\begin{cases}
&2^{k_+}\|P_k(\partial_x(\widetilde{u}_{k_1}\widetilde{v}_{k_2}))\|_{N_k}\lesssim
2^{-\max(|k|,|k_1|,|k_2|)/4} \cdot
2^{{k_1}_+}\|\widetilde{u}_{k_1}\|_{F_{k_1}}\cdot
2^{{k_2}_+}\|\widetilde{v}_{k_2}\|_{F_{k_2}}\\
&\text{if }|k_1-k_2|\leq 4\text{ and }k\leq \min(k_1,k_2)-30.
\end{cases}
\end{equation}
The bound \eqref{on1} follows from \eqref{th1} and \eqref{th2}.

Consider now part (b) of the proposition. We fix extensions
$\widetilde{u}\in C_0(\R:\overline{\E}^\infty)$ of $u$ and
$\widetilde{v}\in C(\R:\E^\infty)$ of $v$ such that
$\|P_k(\widetilde{u})\|_{\overline{F}_k}\leq
2\|P_k(u)\|_{\overline{F}_k(T)}$ and
$\|P_k(\widetilde{v})\|_{F_k}\leq 2\|P_k(v)\|_{F_k(T)}$, $k\in\Z$.
Let $\widetilde{u}_k=P_k(\widetilde{u})$ and
$\widetilde{v}_k=P_k(\widetilde{v})$, $k\in\Z$. It follows from
Lemma \ref{Lemmah1}, Lemma \ref{Lemmah2}, and Lemma \ref{Lemmah4}
that
\begin{equation*}
\begin{split}
(1+2^{-5k/8})\|P_k(\partial_x(\widetilde{u}_{k_1}&
\widetilde{v}_{k_2}))\|_{\overline{N}_k}\lesssim
2^{-\max(|k|,|k_1|,|k_2| )/8}\\
&(1+2^{-5k_1/8})\|\widetilde{u}_{k_1}\|_{\overline{F}_{k_1}}\cdot
(1+2^{k_2})\|v_{k_2}\|_{F_{k_2}}
\end{split}
\end{equation*}
if $|k_1-k|\geq 5$. It follows from Lemma \ref{Lemmah1}, Lemma
\ref{Lemmah2}, and  Lemma \ref{Lemmah3} that
\begin{equation*}
\begin{split}
(1+2^{-5k/8})\|P_k(\partial_x(\widetilde{u}_{k_1}&\widetilde{v}_{k_2}))\|_{\overline{N}_k}\lesssim
2^{-|k_2|/8}(1+2^{-5k_1/8})\|\widetilde{u}_{k_1}\|_{\overline{F}_{k_1}}\cdot
(1+2^{k_2})\|v_{k_2}\|_{F_{k_2}}
\end{split}
\end{equation*}
if $|k_1-k|\leq  4$. The proposition follows.
\end{proof}

The last ingredients in the proof of Theorem \ref{Main1} are
energy estimates. For part (a) we need the following:

\begin{proposition}\label{Lemmad1}
Assume that $T\in(0,1]$ and $u\in C([-T,T]:\E^\infty)$ is a
solution of the initial value problem
\begin{equation}\label{sd8}
\begin{cases}
\partial_tu+\partial_x^3u-\partial_x^{-1}\partial_y^2u+\partial_x(u^2/2)=0\text{ on }\R^2\times[-T,T];\\
u(0)=\phi.
\end{cases}
\end{equation}
Then, for $\sigma\in\{1,2,3\}$ we have
\begin{equation}\label{sd1}
\|u\|^2_{\B^\sigma(T)}\lesssim
\|\phi\|_{\E^\sigma}^2+\|u\|_{\F^1(T)}\cdot
\|u\|_{\F^\sigma(T)}^2.\end{equation}
\end{proposition}

The linearized equation lacks the full set of symmetries of the
nonlinear equation. Consequently, we have good estimates for it
only at the $L^2$ level:

\begin{proposition}\label{Lemmak3}
Assume $T\in(0,1]$, $u\in\overline{\F}^0(T)\cap\F^1(T)$, $v\in
\F^1(T)$ and
\begin{equation}\label{sk8}
\begin{cases}
\partial_tu+\partial_x^3u-\partial_x^{-1}\partial_y^2u+\partial_x(uv)=0\text{ on }\R^2\times(-T,T);\\
u(0)=\phi,
\end{cases}
\end{equation}
Then
\begin{equation}\label{sk9}
\|u\|^2_{\overline{\B}^0(T)}\lesssim
\|\phi\|_{\overline{\E}^0}^2+\|v\|_{\F^1(T)}\cdot
\|u\|_{\overline{\F}^0(T)}^2.
\end{equation}
\end{proposition}

Finally, to estimate differences of solutions in $\F^1$ we need to
differentiate the difference equation. To get bounds for this
equation we need a stronger version of Proposition \ref{Lemmak3}.

\begin{proposition}\label{Lemmak4}
Assume $T\in(0,1]$, $u\in \overline{\F}^0(T)$, $u=P_{\geq
-10}(u)$, $v\in \F^1(T)$,
$w_1,w_2,w_3\in\overline{\F}^0(T)\cap\F^1(T)$, $w'_1,w'_2,w'_3\in
\overline{\F}^0(T)$, $h\in \overline{\F}^0$, $h=P_{\leq -5}(h)$,
and
\begin{equation}\label{su8}
\begin{cases}
\partial_tu+\partial_x^3u-\partial_x^{-1}\partial_y^2u=P_{\geq -10}(v\cdot \partial_xu)+\sum\limits_{m=1}^3P_{\geq -10}(w_m\cdot w'_m)+P_{\geq -10}(h);\\
u(0)=\phi,
\end{cases}
\end{equation}
on $\R^2\times(-T,T)$. Then
\begin{equation}\label{su9}
\begin{split}
\|u\|^2_{\overline{\B}^0(T)}&\lesssim
\|\phi\|_{\overline{\E}^0}^2+\|v\|_{\F^1(T)}\cdot
\|u\|_{\overline{\F}^0(T)}^2+\sum_{m=1}^3\|u\|_{\overline{\F}^0(T)}\|w_m\|_{\overline{\F}^0(T)}\|w'_m\|_{\overline{\F}^0(T)}.
\end{split}
\end{equation}
\end{proposition}

We observe that Proposition \ref{Lemmak3} follows from Proposition
\ref{Lemmak4}: let $u'=P_{\geq -10}u$ and observe that, using
\eqref{sk8} and the definitions
\begin{equation*}
\begin{cases}
&\|u\|_{\overline{\B}^0(T)}^2\lesssim
\|u'\|_{\overline{\B}^0(T)}^2+C\|\phi\|_{\overline{\E}^0}^2;\\
&\partial_x(uv)=v\cdot \partial_xu'+P_{\geq -10}v\cdot
\partial_xP_{\leq -11}u+u\cdot\partial_xv+P_{\leq -11}v\cdot \partial_xP_{\leq
-11}u.
\end{cases}
\end{equation*}
We prove Proposition \ref{Lemmad1} and Proposition \ref{Lemmak4}
in section \ref{energyproof}.

\section{Proof of the main theorem}\label{proof1}

In this section we use the linear, bilinear and energy estimates
in the previous section to prove Theorem \ref{Main1}. Our starting
point is a well-posedness result for more regular solutions:

\begin{proposition}\label{Lemmab7}
Assume $\phi\in \E^\infty$. Then there is $T=T(  \|\phi\|_{\E^3})\in(0,1]$ and a unique solution $u=S^\infty_T(\phi)\in C([-T,T]:\E^\infty)$ of the initial value problem
\begin{equation}\label{er1}
\begin{cases}
\partial_tu+\partial_x^3u-\partial_x^{-1}\partial_y^2u+\partial_x(u^2/2)=0\text{ on }\mathbb{R}^2\times(-T,T);\\
u(0)=\phi.
\end{cases}
\end{equation}
In addition, for any $\sigma\geq  3$
\begin{equation*}
\sup_{t\in[-T,T]}\|u(t)\|_{\E^\sigma}\leq C(\sigma,\|\phi\|_{\E^\sigma},\sup_{t\in[-T,T]}\|u(t)\|_{\E^3}).
\end{equation*}
\end{proposition}

Proposition \ref{Lemmab7} follows by standard energy estimates (see \cite{IoNu}), since
\[
\|\phi\|_{L^\infty}+\|\partial_x\phi\|_{L^\infty}\lesssim\|\phi\|_{\E^3}.
\]

To prove Theorem \ref{Main1} (a), by scaling we may assume that
\begin{equation}\label{er100}
\|\phi\|_{\E^1}\leq \eps_0\ll1.
\end{equation}
The uniqueness part of Theorem
\ref{Main1} (a) follows from Proposition \ref{Lemmab7}. For global
existence, in view of the conservation laws \eqref{conserve}, we only
need to construct the solution on the time interval $[-1,1]$.
In view of Proposition \ref{Lemmab7}, it suffices to prove that if
$T\in(0,1]$ and $u\in C([-T,T]:\E^\infty)$ is a solution of
\eqref{er1} with  $\|\phi\|_{\E^1}\leq \eps_0\ll1$ then
\begin{equation}\label{er2}
\sup_{t\in[-T,T]}\|u(t)\|_{\E^3}\lesssim \|\phi\|_{\E^3}.
\end{equation}

We first use a continuity argument to establish an $\F^1$ bound on $u$
in the interval $[-T,T]$.  By taking $\sigma=1$, it follows from
Proposition \ref{Lemmaa3} (a), Proposition \ref{Lemmab1} (a), and Proposition
\ref{Lemmad1} that for any $T'\in[0,T]$ we have
\begin{equation}\label{er3}
\begin{cases}
&\|u\|_{\F^1(T')}\lesssim \|u\|_{\B^1(T')}+\|\partial_x(u^2)\|_{\N^1(T')};\\
&\|\partial_x(u^2)\|_{\N^1(T')}\lesssim \|u\|_{\F^1(T')}^2;\\
&\|u\|_{\B^1(T')}^2\lesssim \|\phi\|_{\E^1}^2+\|u\|_{\F^1(T')}^3.
\end{cases}
\end{equation}
We denote $X(T')=\|u\|_{\B^1(T')}+\|\partial_x(u^2)\|_{\N^1(T')}$ and
eliminate $\|u\|_{\F^1(T')}$ to obtain
\begin{equation}\label{er4}
X(T')^2\lesssim \|\phi\|_{\E^1}^2+X(T')^3+X(T')^4.
\end{equation}
Assuming that $X(T')$ is continuous and satisfies
\begin{equation} \label{xtcont}
\lim_{T'\to 0} X(T') \lesssim \|\phi\|_{\E^1}
\end{equation}
if $\eps_0$ is sufficiently small, we would conclude
that $X(T') \lesssim \|\phi\|_{\E^1}$. Using \eqref{er3},
\begin{equation}\label{er6}
\|u\|_{\F^1(T)}\lesssim\|\phi\|_{\E^1},
\end{equation}

To obtain \eqref{xtcont} and the continuity of $X(T')$ we first
observe that for $u\in C([-T,T]:\E^\infty)$ the mapping $T'\to
\|u\|_{\B^1(T')}$ is increasing and continuous on the interval
$[-T,T]$ and
\begin{equation*}
\lim_{T'\to 0}\|u\|_{\B^1(T')}\lesssim \|\phi\|_{\E^1}.
\end{equation*}
The similar properties of $\|\partial_x(u^2)\|_{\N^1(T')}$ are
obtained by applying the following lemma to $v =  \partial_x(u^2)$.

\begin{lemma}\label{Lemmar1}
Assume $T\in(0,1]$ and $v\in C([-T,T]:\E^\infty)$. Then the mapping
$T'\to\|v\|_{\N^1(T')}$ is increasing and continuous on the interval
$[0,T]$ and
\begin{equation}\label{er30}
\lim_{T'\to 0}\|v\|_{\N^1(T')}=0.
\end{equation}
\end{lemma}

\begin{proof}[Proof of Lemma \ref{Lemmar1}] In view of the definitions and \eqref{sb21},
\begin{equation}\label{er31}
\|v\|_{\N^1(T')}\lesssim \|p(\xi,\mu)\cdot \mathcal{F}(v\cdot \mathbf{1}_{[-T',T']}(t))
\|_{L^2}\lesssim T'^{1/2}\sup_{t\in[-T',T']}\|v(t)\|_{\E^1}.
\end{equation}
The limit in \eqref{er30} follows. It remains to prove the
continuity of the mapping $T'\to\|v\|_{\N^1(T')}$ at some point
$T'_0\in(0,T]$. We fix $\eps>0$. Let $D_r(v)(x,y,t)=v(x,y,t/r)$,
$r\in[1/2,2]$. Using \eqref{er31} we have
\begin{equation*}
\|v\|_{\N^1(T')}-\|D_{T'/T'_0}(v)\|_{\N^1(T')}\lesssim \sup_{t\in[-T',T']}\|v(t)-D_{T'/T'_0}(v)(t)\|_{\E^1}\leq \eps,
\end{equation*}
for any $T'\in(0,T]$ sufficiently close to $T'_0$. Also, using the
definitions
\begin{equation*}
\lim_{r\to 1}\|D_r(v)\|_{\N^1(rT'_0)}=\|v\|_{\N^1(T'_0)},
\end{equation*}
which completes the proof of the lemma.
\end{proof}

To prove \eqref{er2} we combine again Proposition \ref{Lemmaa3} (a),
Proposition \ref{Lemmab1} (a), and Proposition \ref{Lemmad1} with
$\sigma=2,3$ to conclude that
\begin{equation}\label{er8}
\begin{cases}
&\|u\|_{\F^\sigma(T)}\lesssim \|u\|_{\B^\sigma(T)}+\|\partial_x(u^2)\|_{\N^\sigma(T)};\\
&\|\partial_x(u^2)\|_{\N^\sigma(T)}\lesssim \|u\|_{\F^1(T)}\cdot \|u\|_{\F^\sigma(T)};\\
&\|u\|_{\B^\sigma(T)}^2\lesssim \|\phi\|_{\E^\sigma}^2+\|u\|_{\F^1(T)}\cdot \|u\|_{\F^\sigma}^2.
\end{cases}
\end{equation}
Using \eqref{er100} and \eqref{er6}, it follows that
\begin{equation} \label{e2bd}
\|u\|_{\F^\sigma(T)}\lesssim 
\|\phi\|_{\E^\sigma}\text{ for }\sigma\in\{2,3\}.
\end{equation}
Then the inequality \eqref{er2} follows from Lemma \ref{Lemmaa2}.

We prove now Theorem \ref{Main1} (b). Assume $\phi\in
\mathbf{E}^1$ is fixed,
\begin{equation*}
\{\phi_n:n\in\Z_+\}\subseteq \E^\infty\text{ and }\lim_{n\to\infty}\phi_n=\phi\text{ in }\E^1.
\end{equation*}
It suffices to prove that the sequence $S^\infty_T(\phi_n)\in C([-T,T]:\E^\infty)$ is a Cauchy sequence in $C([-T,T]:\E^1)$. By scaling, we may assume
\begin{equation}\label{amy0}
\|\phi\|_{\E^1}\leq \eps_0\ll1\text{ and }\|\phi_n\|_{\E^1}\leq\eps_0\ll 1\text{ for any }n\in\Z_+.
\end{equation}
Using the conservation laws \eqref{conserve} it suffices to prove that for any $\delta>0$ there is $M_\delta$ such that
\begin{equation}\label{amy1}
\sup_{t\in[-1,1]}\|S^\infty(\phi_m)(t)-S^\infty(\phi_n)(t)\|_{\E^1}\leq \delta\text{ for any }m,n\geq M_\delta.
\end{equation}
For $K\in\Z_+$ let
\begin{equation*}
\phi_n^K=P_{\leq  K}\phi_n.
\end{equation*}

We show first that for any $K\in\Z_+$ there is $M_{\delta,K}$ such that
\begin{equation}\label{amy2}
\sup_{t\in[-1,1]}\|S^\infty(\phi_m^K)(t)-S^\infty(\phi_n^K)(t)\|_{\E^1}\leq \delta\text{ for any }m,n\geq M_{\delta,K}.
\end{equation}
Using Theorem \ref{Main1} (a)
\begin{equation*}
\sup_{t\in[-1,1]}\|S^\infty(\phi_n^K)\|_{\E^{10}}\leq C(K)\text{ for any }n\in \Z_+.
\end{equation*}
Standard energy estimates for the difference equation show that
\begin{equation*}
\sup_{t\in[-1,1]}\|S^\infty(\phi_m^K)(t)-S^\infty(\phi_n^K)(t)\|_{\E^1}\leq C'(K)\cdot \|\phi_m-\phi_n\|_{\E^1},
\end{equation*}
and \eqref{amy2} follows.

We show now that for any $\delta>0$ there are $K\in\Z_+$ and $M_\delta$ such that
\begin{equation}\label{amy3}
\sup_{t\in[-1,1]}\|S^\infty(\phi_n)(t)-S^\infty(\phi_n^K)(t)\|_{\E^1}\leq \delta\text{ for any }n\geq M_{\delta}.
\end{equation}
The main bound \eqref{amy1} would follow from \eqref{amy2} and \eqref{amy3}. To prove \eqref{amy3} we need to estimate differences of solutions. We summarize our main result
in Proposition \ref{diffest} below. The bound \eqref{amy3} follows from \eqref{vf1} and Lemma \ref{Lemmaa2}:
\begin{equation*}
\begin{split}
\sup_{t\in[-1,1]}\|S^\infty(\phi_n)(t)-S^\infty(\phi_n^K)(t)\|_{\E^1}&\lesssim \|S^\infty(\phi_n)(t)-S^\infty(\phi_n^K)(t)\|_{\F^1(1)}\\
&\lesssim\|\phi_n-\phi_n^K\|_{\E^1}+C\|\phi_n^K\|_{\E^2}\|\phi_n-\phi_n^K\|_{\overline{\E}^0}\\
&\lesssim\|\phi-\phi_n\|_{\E^1}+\|\phi-P_{\leq K}\phi\|_{\E^1}.
\end{split}
\end{equation*}

\begin{proposition} \label{diffest}
Let $u_1,u_2\in \F^1(1)$ be solutions to \eqref{eq-1} with initial data
$\phi_1, \phi_2 \in \E^\infty$ satisfying
\[
\|\phi_1\|_{\E^1}+\|\phi_2\|_{\E^1}\leq \eps_0\ll 1, \qquad \phi_1-\phi_2 \in \overline{\E}^0.
\]
Then
\begin{equation}\label{vfo}
 \| u_1 - u_2\|_{\overline\F^0(1)}
\lesssim \| \phi_1 - \phi_2\|_{\overline\E^0},
\end{equation}
and
\begin{equation}\label{vf1}
  \| u_1 - u_2\|_{\F^1(1)}
  \lesssim \|\phi_1 -\phi_2 \|_{\E^1} +
  \|\phi_1\|_{\E^2}  \| \phi_1 - \phi_2\|_{\overline\E^0}
\end{equation}
\end{proposition}

\begin{proof}[Proof of Proposition \ref{diffest}]
The difference $v=u_2-u_1$ solves the equation
\begin{equation}\label{amy4}
\begin{cases}
  &\partial_t v+\partial_x^3 v-\partial_x^{-1}\partial_y^2v=
-\partial_x[v(u_1+u_2)/2];\\
  &v(0)=\phi = \phi_2-\phi_1.
\end{cases}
\end{equation}
By \eqref{er6} we can bound
\begin{equation}\label{amy5}
\|u_1\|_{\F^1(1)}+\|u_2\|_{\F^1(1)}\lesssim\eps_0.
\end{equation}
By Proposition \ref{Lemmaa3} (b)\footnote{Clearly,
$\partial_x[v(u_1+u_2)/2]\in \overline{\N}^0(1)$ since
$\partial_x[v(u_1+u_2)/2]\in C([-1,1]:\overline{\E}^0)$.},
Proposition \ref{Lemmab1} (b), and Proposition \ref{Lemmak3}
we obtain
\begin{equation*}
\begin{cases}
&\|v\|_{\overline{\F}^0(1)}\lesssim \|v\|_{\overline{\B}^0(1)}+
\|\partial_x[v(u_1+u_2)/2]\|_{\overline{\N}^0(1)};\\
&\|\partial_x[v(u_1+u_2)/2]\|_{\overline{\N}^0(1)}\lesssim \|v\|_{\overline{\F}^0(1)}( \|u_1\|_{\F^1(1)}+\|u_{2}\|_{\F^1(1)});\\
&\|v\|^2_{\overline{\B}^0(1)}\lesssim \|\phi\|_{\overline{\E}^0}^2+\|v\|^2_{\overline{\F}^0(1)}( \|u_1\|_{\F^1(1)}+\|u_{2}\|_{\F^1(1)}).
\end{cases}
\end{equation*}
Combining this with  \eqref{amy5} we obtain the estimate \eqref{vfo}.

To prove \eqref{vf1} we first use Proposition~\ref{Lemmaa3} (a) and
Proposition \ref{Lemmab1} (a) to obtain
\begin{equation*}
\begin{cases}
&\|v\|_{\F^1(1)}\lesssim \|v\|_{\B^1(1)}+\|\partial_x[v(u_1+u_2)/2]\|_{\N^1(1)};\\
&\|\partial_x[v(u_1+u_2)/2]\|_{\N^1(1)}\lesssim \|v\|_{\F^1(1)}\cdot ( \|u_1\|_{\F^1(1)}+\|u_2\|_{\F^1(1)}).
\end{cases}
\end{equation*}
Since $\|P_{\leq 0} (v)\|_{\B^1(1)}=\|P_{\leq 0} (\phi)\|_{\E^1}$, it follows from \eqref{amy5} that
\begin{equation}\label{amy6}
\|v\|_{\F^1(1)}\lesssim \|P_{\geq 1} (v)\|_{\B^1(1)} + \|\phi\|_{\E^1}
\end{equation}
By \eqref{vfo} and \eqref{e2bd}, for \eqref{vf1} it remains to
prove the estimate
\begin{equation} \label{p1v}
\| P_{\geq 1} (v)\|_{\B^1(1)}^2 \lesssim \|\phi \|_{\E^1}^2
+  \|u_1\|_{\F^2(1)}\cdot \| v \|_{\overline\F^0(1)}(\|P_{\geq 1} (v)\|_{\B^1(1)} + \|\phi\|_{\E^1}).
\end{equation}
In view of the definitions,
\begin{equation}\label{amy30}
  \|P_{\geq 1} (v)\|_{\B^1(1)}\approx \|P_{\geq 1} (\partial_x
  v)\|_{\overline{\B}^0(1)}+
\|P_{\geq 1} (\partial_x^{-1}\partial_y v)\|_{\overline{\B}^0(1)}.
\end{equation}

Using \eqref{amy4} we write the equation for $U_1=P_{\geq -10}(\partial_xv)$ in the form
\begin{equation*}
\begin{cases}
&\partial_tU_1+\partial_x^3U_1-\partial_x^{-1}\partial_y^2U_1=P_{\geq -10}(-u_2\cdot \partial_xU_1)+P_{\geq -10}(G_1);\\
&U_1(0)=P_{\geq -10}(\partial_x\phi),
\end{cases}
\end{equation*}
where
\begin{equation*}
\begin{split}
G_1=&-P_{\geq -10}(u_2)\cdot \partial_x^2P_{\leq -11}(v)-P_{\leq -11}(u_2)\cdot \partial_x^2P_{\leq -11}(v)\\
&-\partial_xv\cdot \partial_x(u_1+u_2)-v\cdot \partial^2_xu_1.
\end{split}
\end{equation*}
It follows from Proposition \ref{Lemmak4} (with $h=-P_{\leq -11}(u_2)\cdot \partial_x^2P_{\leq -11}(v)$) that
\begin{equation*}
\begin{split}
\|P_{\geq -10}&(\partial_xv)\|^2_{\overline{\B}^0(1)}\lesssim \|\partial_x\phi\|_{\overline{\E}^0}^2+\|u_2\|_{\F^1(1)}\|P_{\geq -10}(\partial_xv)\|^2_{\overline{\F}^0(1)}\\
&+\|P_{\geq -10}(\partial_xv)\|_{\overline{\F}^0(1)}\cdot \|\partial_xv\|_{\overline{\F}^0(1)}\cdot \left[ \|\partial_xu_1\|_{\overline{\F}^0(1)}+\|\partial_xu_2\|_{\overline{\F}^0(1)}\right]\\
&+\|P_{\geq -10}(\partial_xv)\|_{\overline{\F}^0(1)}\cdot \|v\|_{\overline{\F}^0(1)}\cdot \|\partial^2_xu_1\|_{\overline{\F}^0(1)}.
\end{split}
\end{equation*}
Using \eqref{amy5} and \eqref{amy6}, it follows that
\begin{equation}\label{amy50}
\begin{split}
\|P_{\geq -10}(\partial_xv)&\|^2_{\overline{\B}^0(1)}\lesssim\|\phi\|_{{\E}^1}^2+\eps_0(\|P_{\geq 1}(v)\|_{\B^1(1)}^2+\|\phi\|_{\E^1}^2)\\
&+( \|P_{\geq 1}(v)\|_{\B^1(1)}+\|\phi\|_{\E^1})\cdot \|v\|_{\overline{\F}^0(1)}\cdot \|u_1\|_{\F^2(1)}.
\end{split}
\end{equation}

Using \eqref{amy4} we write the equation for $U_2=P_{\geq -10} (\partial_x^{-1}\partial_yv)$ in the form
\begin{equation*}
\begin{cases}
&\partial_tU_2+\partial_x^3U_2-\partial_x^{-1}\partial_y^2U_2=P_{\geq -10}(-u_2\cdot \partial_xU_2)+P_{\geq -10}G_2;\\
&U_2(0)=P_{\geq -10}(\partial_x^{-1}\partial_y\phi),
\end{cases}
\end{equation*}
where
\begin{equation*}
\begin{split}
G_2=&-P_{\geq -10}(u_2)\cdot \partial_xP_{\leq -11}(\partial_x^{-1}\partial_yv)-P_{\leq  -11}(u_2)\cdot \partial_xP_{\leq -11}(\partial_x^{-1}\partial_yv)-v\cdot \partial_yu_1.
\end{split}
\end{equation*}
It follows from Proposition \ref{Lemmak4} (with $h=-P_{\leq  -11}(u_2)\cdot \partial_xP_{\leq -11}(\partial_x^{-1}\partial_yv)$) that
\begin{equation*}
\begin{split}
&\|P_{\geq
  -10}(\partial_x^{-1}\partial_yv)\|^2_{\overline{\B}^0(1)}\lesssim \|P_{\geq -10}(\partial_x^{-1}\partial_y\phi)\|_{\overline{\E}^0}^2+\|u_2\|_{\F^1(1)}\|v\|^2_{\F^1(1)}\\
&+\|P_{\geq -10}(\partial_x^{-1}\partial_yv)\|_{\overline{\F}^0(1)}\cdot \|v\|_{\overline{\F}^0(1)}\cdot \|\partial_x(\partial_x^{-1}\partial_yu_1)\|_{\overline{\F}^0(1)}.
\end{split}
\end{equation*}
Using \eqref{amy5} and \eqref{amy6}, it follows that
\begin{equation}\label{amy51}
\begin{split}
\|P_{\geq -10}(\partial_x^{-1}&\partial_yv)\|^2_{\overline{\B}^0(1)}\lesssim \|\phi\|_{\E^1}^2+\eps_0(\|P_{\geq 1}(v)\|_{\B^1(1)}^2+\|\phi\|_{\E^1}^2)\\
&+( \|P_{\geq 1}(v)\|_{\B^1(1)}+\|\phi\|_{\E^1})\cdot \|v\|_{\overline{\F}^0(1)}\cdot \|u_1\|_{{\F}^2(1)}.
\end{split}
\end{equation}
We add up \eqref{amy50} and \eqref{amy51} and use \eqref{amy30}. The bound \eqref{p1v} follows.
\end{proof}

\section{$L^2$ bilinear estimates}\label{L2bi}

For $k\in\Z$ and $l,j\in\R$ let
\begin{equation*}
D_{k,l,j}=\{(\xi,\mu,\tau):\,\xi\in \widetilde{I}_k,\,|\mu|\leq
2^l,\,|\tau-\omega(\xi,\mu)|\leq 2^j\},
\end{equation*}
and let $D_{k,\infty,j}=\cup_{l\in\Z}D_{k,l,j}$.

\newtheorem{Main9}{Lemma}[section]
\begin{Main9}\label{Main9}
(a) Assume $k_1,k_2,k_3\in \Z$, $j_1,j_2,j_3\in\Z_+$, and
$f_i:\R^3\to\R_+$ are $L^2$ functions supported in
$D_{k_i,\infty,j_i}$, $i=1,2,3$. If
\begin{equation}\label{jj1.1}
\max(j_1,j_2,j_3)\leq k_1+k_2+k_3-20
\end{equation}
then
\begin{equation}\label{jj1.6}
\int_{\R^3}(f_1\ast f_2)\cdot f_3\lesssim   2^{(j_1+j_2+j_3)/2}\cdot
2^{-(k_1+k_2+k_3)/2}\cdot \|f_1\|_{L^2}\|f_2\|_{L^2}\|f_3\|_{L^2}.
\end{equation}

(b)Assume $k_1,k_2,k_3\in \Z$, $l_1,l_2,l_3\in\Z$,
$j_1,j_2,j_3\in\Z_+$, and $f_i:\R^3\to\R_+$ are $L^2$ functions
supported in $D_{k_i,l_i,j_i}$, $i=1,2,3$. Then
\begin{equation}\label{jj2.22}
\int_{\R^3}(f_1\ast f_2)\cdot f_3\lesssim 
 2^{[\min(k_1,k_2,k_3)+\min(l_1,l_2,l_3)+\min(j_1,j_2,j_3)]/2}\cdot
\|f_1\|_{L^2}\|f_2\|_{L^2}\|f_3\|_{L^2}.
\end{equation}
\end{Main9}

\begin{proof}[Proof of Lemma \ref{Main9}]
Part (b) follows easily from the Minkowski inequality. Part (a) is proved also in \cite{CoIoKeSt}; we reproduce the proof here for the sake of completeness. We observe that
\begin{equation}\label{mola6}
\int_{\R^3}(f_1\ast f_2)\cdot f_3=\int_{\R^3}(\widetilde{f}_1\ast
f_3)\cdot f_2=\int_{\R^3}(\widetilde{f}_2\ast f_3)\cdot f_1,
\end{equation}
where $\widetilde{f}_i(\xi,\mu,\tau)=f_i(-\xi,-\mu,-\tau)$,
$i=1,2$. In view of the symmetry of \eqref{jj1.6} we may assume
\begin{equation}\label{jj0}
j_3=\max(j_1,j_2,j_3).
\end{equation}
We define
$f_i^\#(\xi,\mu,\theta)=f_i(\xi,\mu,\theta+\omega(\xi,\mu))$,
$i=1,2,3$, $\|f_i^\#\|_{L^2}=\|f_i\|_{L^2}$. We rewrite the
left-hand side of \eqref{jj1.6} in the form
\begin{equation}\label{jj3}
\begin{split}
&\int_{\R^6}f_1^\#(\xi_1,\mu_1,\theta_1)\cdot f_2^\#(\xi_2,\mu_2,\theta_2)\\
&\times
f_3^\#(\xi_1+\xi_2,\mu_1+\mu_2,\theta_1+\theta_2+\Omega((\xi_1,\mu_1),(\xi_2,\mu_2)))\,d\xi_1d\xi_2d\mu_1d\mu_2d\theta_1d\theta_2,
\end{split}
\end{equation}
where
\begin{equation}\label{jj2.6}
\begin{split}
\Omega((\xi_1,\mu_1),(\xi_2,\mu_2))&=-\omega(\xi_1+\xi_2,\mu_1+\mu_2)+\omega(\xi_1,\mu_1)+\omega(\xi_2,\mu_2)\\
&=\frac{-\xi_1\xi_2}{\xi_1+\xi_2}\Big[(\sqrt{3}\xi_1+\sqrt{3}\xi_2)^2-\Big(\frac{\mu_1}{\xi_1}-\frac{\mu_2}{\xi_2}\Big)^2\Big].
\end{split}
\end{equation}
The functions $f_i^\#$ are supported in the sets
$\{\xi,\mu,\theta):\,\xi\in
\widetilde{I}_{k_i},\,\mu\in\R,\,|\theta|\leq 2^{j_i}\}$.

We will prove that if  $g_i:\R^2\to\R_+$ are $L^2$ functions
supported in $\widetilde{I}_{k_i}\times\R$, $i=1,2$, and
$g:\R^3\to\R_+$ is an $L^2$ function supported in
$\widetilde{I}_k\times\R\times [-2^j,2^j]$, $j\leq k_1+k_2+k-15$,
then
\begin{equation}\label{jj5}
\begin{split}
\int_{\R^4}g_1(\xi_1,\mu_1)\cdot &g_2(\xi_2,\mu_2)\cdot g(\xi_1+\xi_2,\mu_1+\mu_2,\Omega((\xi_1,\mu_1),(\xi_2,\mu_2)))\,d\xi_1d\xi_2d\mu_1d\mu_2\\
&\lesssim   2^{j/2}\cdot
2^{-(k_1+k_2+k)/2}\cdot\|g_1\|_{L^2}\|g_2\|_{L^2}\|g\|_{L^2}.
\end{split}
\end{equation}
This suffices for \eqref{jj1.6}, in view of \eqref{jj0} and
\eqref{jj3}.

To prove \eqref{jj5}, we observe\footnote{There are four identical
integrals of this type.} first that we may assume that the
integral in the left-hand side of \eqref{jj5} is taken over the
set
\begin{equation*}
\mathcal{R}_{++}=\{(\xi_1,\mu_1,\xi_2,\mu_2):\,\xi_1+\xi_2\geq
0\text{ and }\mu_1/ \xi_1-\mu_2/ \xi_2\geq 0\}.
\end{equation*}
Using the restriction $j\leq k_1+k_2+k-15$ and \eqref{jj2.6}, we
may assume also that the integral in the left-hand side of
\eqref{jj5} is taken over the set
\begin{equation*}
\widetilde{\mathcal{R}}_{++}=\{(\xi_1,\mu_1,\xi_2,\mu_2)\in\mathcal{R}_{++}:|\sqrt{3}(\xi_1+\xi_2)|-|\mu_1/
\xi_1-\mu_2/ \xi_2|\leq 2^{-10}|\xi_1+\xi_2|\}.
\end{equation*}
To summarize, it suffices to prove that
\begin{equation}\label{jj11}
\begin{split}
\int_{\widetilde{\mathcal{R}}_{++}}g_1(\xi_1,\mu_1)\cdot
&g_2(\xi_2,\mu_2)\cdot
g(\xi_1+\xi_2,\mu_1+\mu_2,\Omega((\xi_1,\mu_1),(\xi_2,\mu_2)))\,d\xi_1d\xi_2d\mu_1d\mu_2\\
&\lesssim   2^{j/2}\cdot
2^{-(k_1+k_2+k)/2}\cdot\|g_1\|_{L^2}\|g_2\|_{L^2}\|g\|_{L^2}.
\end{split}
\end{equation}

We make the changes of variables
\begin{equation*}
\mu_1=\sqrt{3}\xi_1^2+\beta_1\xi_1\text{ and
}\mu_2=-\sqrt{3}\xi_2^2+\beta_2\xi_2,
\end{equation*}
with $d\mu_1d\mu_2=\xi_1\xi_2\,d\beta_1d\beta_2$. The left-hand
side of \eqref{jj11} is bounded by
\begin{equation}\label{jj12}
\begin{split}
& C2^{k_1+k_2}\int_{S}g_1(\xi_1,\sqrt{3}\xi_1^2+\beta_1\xi_1)\cdot g_2(\xi_2,-\sqrt{3}\xi_2^2+\beta_2\xi_2)\\
&\times
g(\xi_1+\xi_2,\sqrt{3}\xi_1^2-\sqrt{3}\xi_2^2+\beta_1\xi_1+\beta_2\xi_2,\widetilde{\Omega}((\xi_1,\beta_1),(\xi_2,\beta_2)))\,d\xi_1d\xi_2d\beta_1d\beta_2,
\end{split}
\end{equation}
where
\begin{equation}\label{jj10}
S=\{(\xi_1,\beta_1,\xi_2,\beta_2):\xi_1+\xi_2\geq 0\text{ and
}|\beta_1-\beta_2|\leq 2^{-10}(\xi_1+\xi_2)\},
\end{equation}
and
\begin{equation}\label{jj15}
\widetilde{\Omega}((\xi_1,\beta_1),(\xi_2,\beta_2))=\xi_1\xi_2
(\beta_1-\beta_2)\Big(2\sqrt{3}+\frac{\beta_1-\beta_2}{\xi_1+\xi_2}\Big).
\end{equation}

We define the functions $h_i:\R^2\to\R_+$ supported in
$\widetilde{I}_{k_i}\times\R$, $i=1,2$,
\begin{equation*}
\begin{cases}
&h_1(\xi_1,\beta_1)=2^{k_1/2}\cdot g_1(\xi_1,\sqrt{3}\xi_1^2+\beta_1\xi_1);\\
&h_2(\xi_2,\beta_2)=2^{k_2/2}\cdot
g_2(\xi_2,-\sqrt{3}\xi_2^2+\beta_2\xi_2),
\end{cases}
\end{equation*}
with $\|h_i\|_{L^2}\approx\|g_i\|_{L^2}$. Thus, for \eqref{jj5} it
suffices to prove that
\begin{equation}\label{jj20}
\begin{split}
&2^{(k_1+k_2)/2}\int_{S}h_1(\xi_1,\beta_1)\cdot h_2(\xi_2,\beta_2)\\
&\times g(\xi_1+\xi_2,\sqrt{3}\xi_1^2-\sqrt{3}\xi_2^2+\beta_1\xi_1+\beta_2\xi_2,\widetilde{\Omega}((\xi_1,\beta_1),(\xi_2,\beta_2)))\,d\xi_1d\xi_2d\beta_1d\beta_2\\
&\lesssim   2^{j/2}\cdot
2^{-(k_1+k_2+k)/2}\cdot\|h_1\|_{L^2}\|h_2\|_{L^2}\|g\|_{L^2}.
\end{split}
\end{equation}

To prove \eqref{jj20}, we may assume without loss of generality
that
\begin{equation}\label{jj11.1}
k_1\leq k_2.
\end{equation}
We make the change of variables $\beta_1=\beta_2+\beta$. In view
of \eqref{jj10}, \eqref{jj15}, and the restriction on the support
of $g$, we may assume $|\beta|\leq 2^{j-k_1-k_2+4}$. Thus, the
integral in the left-hand side of \eqref{jj20} is equal to
\begin{equation}\label{jj21}
\begin{split}
&2^{(k_1+k_2)/2}\int_{\widetilde{S}}h_1(\xi_1,\beta+\beta_2)\cdot h_2(\xi_2,\beta_2)\cdot\mathbf{1}_{[-1,1]}( \beta/2^{j-k_1-k_2+4} )\\
&\times
g(\xi_1+\xi_2,A(\xi_1,\xi_2,\beta)+\beta_2(\xi_1+\xi_2),B(\xi_1,\xi_2,\beta))\,d\xi_1d\xi_2d\beta
d\beta_2,
\end{split}
\end{equation}
where
$\widetilde{S}=\{(\xi_1,\xi_2,\beta,\beta_2)\in\R^4:\xi_1+\xi_2\geq
0\text{ and }|\beta|\leq 2^{-10}(\xi_1+\xi_2)\}$, and
\begin{equation}\label{jj25}
\begin{cases}
&A(\xi_1,\xi_2,\beta)=\sqrt{3}\xi_1^2-\sqrt{3}\xi_2^2+\beta \xi_1;\\
&B(\xi_1,\xi_2,\beta)=\xi_1\xi_2\beta\cdot
(2\sqrt{3}+\beta/(\xi_1+\xi_2)).
\end{cases}
\end{equation}
Let $j'=j-k_1-k_2+4$ and decompose, for $i=1,2$,
\begin{equation*}
h_i(\xi',\beta')=\sum_{m\in\Z}h_i(\xi',\beta')\cdot
\mathbf{1}_{[0,1)}(\beta'/2^{j'}-m)=\sum_{m\in\Z}h_i^m(\xi',\beta').
\end{equation*}
The expression in \eqref{jj21} is dominated by
\begin{equation}\label{jj26}
\begin{split}
&C2^{(k_1+k_2)/2}\sum_{|m-m'|\leq 4}\int_{\widetilde{S}}h^m_1(\xi_1,\beta+\beta_2)\cdot h^{m'}_2(\xi_2,\beta_2)\\
&\times
g(\xi_1+\xi_2,A(\xi_1,\xi_2,\beta)+\beta_2(\xi_1+\xi_2),B(\xi_1,\xi_2,\beta))\,d\xi_1d\xi_2d\beta
d\beta_2.
\end{split}
\end{equation}
Also, for $i=1,2$,
\begin{equation*}
\|h_i\|_{L^2}=\big[\sum_{m\in\Z}\|h_i^m\|_{L^2}^2\big].
\end{equation*}
Thus, to prove \eqref{jj20}, we may assume $h_1=h_1^m$ and
$h_2=h_2^{m'}$ for some fixed $m,m'\in\Z$ with $|m-m'|\leq 4$. To
summarize, it suffices to prove that if $F_i:\R^2\to[0,\infty)$
are $L^2$ functions supported in $\widetilde{I}_{k_i}\times\R$,
$g$ is as before, and $m\in\Z$ then
\begin{equation}\label{jj30}
\begin{split}
&2^{(k_1+k_2)/2}\int_{\widetilde{S}}F_1(\xi_1,\beta+\beta_2)\cdot F_2(\xi_2,\beta_2)\cdot \mathbf{1}_{[m-1,m+1]}(\beta_2/2^{j'})\\
&\times g(\xi_1+\xi_2,A(\xi_1,\xi_2,\beta)+\beta_2(\xi_1+\xi_2),B(\xi_1,\xi_2,\beta))\,d\xi_1d\xi_2d\beta d\beta_2\\
&\lesssim   2^{j/2}\cdot
2^{-(k_1+k_2+k)/2}\cdot\|F_1\|_{L^2}\|F_2\|_{L^2}\|g\|_{L^2}.
\end{split}
\end{equation}
To prove  \eqref{jj30} we use the Minkowski inequality in the
variables $(\xi_1,\xi_2,\beta)$: with
\begin{equation*}
S'=\{(\xi_1,\xi_2,\beta)\in\R^3:\xi_i\in
\widetilde{I}_{k_i},\,\xi_1+\xi_2\geq 0,\,|\beta|\leq
2^{-10}(\xi_1+\xi_2)\},
\end{equation*}
the left-hand side of \eqref{jj30} is dominated by
\begin{equation}\label{jj35}
\begin{split}
C2^{(k_1+k_2)/2}\int_\R\mathbf{1}_{[m-1,m+1]}(\beta_2/2^{j'})\cdot \Big(\int_{S'}|F_1(\xi_1,\beta+\beta_2)\cdot F_2(\xi_2,\beta_2)|^2\,d\xi_1d\xi_2d\beta\Big)^{1/2}\\
\times\Big(\int_{S'}|g(\xi_1+\xi_2,A(\xi_1,\xi_2,\beta)+\beta_2(\xi_1+\xi_2),B(\xi_1,\xi_2,\beta))|^2\,d\xi_1d\xi_2d\beta
\Big)^{1/2}\,d\beta_2.
\end{split}
\end{equation}
For \eqref{jj30}, it is easy to see that it suffices to prove that
\begin{equation}\label{jj40}
\begin{split}
\Big(\int_{S'}|g(\xi_1+\xi_2,A(\xi_1,\xi_2,\beta)+\beta_2(\xi_1+\xi_2),B(\xi_1,\xi_2,\beta))|^2\,d\xi_1d\xi_2d\beta \Big)^{1/2}\\
\lesssim   2^{-(k_1+k_2+k)/2}\|g\|_{L^2}.
\end{split}
\end{equation}
for any $\beta_2\in\R$. Indeed, assuming \eqref{jj40}, we can
bound the expression in \eqref{jj35} by
\begin{equation*}
C2^{(k_1+k_2)/2}\int_\R\mathbf{1}_{[m-1,m+1]}(\beta_2/2^{j'})\cdot
\|F_1\|_{L^2}\|F_2(.,\beta_2)\|_{L^2_{\xi_2}}\cdot
2^{-(k_1+k_2+k)/2}\|g\|_{L^2}\,d\beta_2,
\end{equation*}
which suffices since $2^{j'/2}2^{(k_1+k_2)/2}\approx 2^{j/2}$.

Finally, to prove \eqref{jj40}, we may assume first that
$\beta_2=0$. We examine \eqref{jj25} and make the change of
variable $\beta=\sqrt{3}(\xi_1+\xi_2)\cdot \nu$. The left-hand
side of \eqref{jj40} is dominated by
\begin{equation}\label{jj41}
C\Big(2^{k}\int_{S''}|g(\xi_1+\xi_2,\sqrt{3}(\xi_1+\xi_2)(\xi_1-\xi_2+\nu\xi_1),3\xi_1\xi_2(\xi_1+\xi_2)\nu(2+\nu))|^2\,d\xi_1d\xi_2d\nu
\Big)^{1/2},
\end{equation}
where $S''=\{(\xi_1,\xi_2,\nu)\in\R^3:\xi_i\in
\widetilde{I}_{k_i},\,|\nu|\leq 2^{-10}\}$. We define the function
\begin{equation*}
h(\xi,x,y)=2^{2k}\cdot |g(\xi,\sqrt{3}\xi\cdot x,3\xi\cdot y)|^2,
\end{equation*}
so $\|h\||_{L^1}\approx \|g\|_{L^2}^2$. The expression in
\eqref{jj41} is dominated by
\begin{equation*}
C2^{-k/2}\Big(\int_{S''}|h(\xi_1+\xi_2,\xi_1-\xi_2+\nu\xi_1,\xi_1\xi_2\cdot
\nu(2+\nu))|\,d\xi_1d\xi_2d\nu \Big)^{1/2}.
\end{equation*}
Therefore, it remains to prove that
\begin{equation*}
\int_{S''}|h(\xi_1+\xi_2,\xi_1-\xi_2+\nu\xi_1,\xi_1\xi_2\cdot
\nu(2+\nu))|\,d\xi_1d\xi_2d\nu\lesssim   2^{-(k_1+k_2)}\|h\|_{L^1}
\end{equation*}
for any function $h\in L^1(\R^3)$. This is clear since the
absolute value of the determinant of the change of variables
$(\xi_1,\xi_2,\nu)\to
[\xi_1+\xi_2,\xi_1-\xi_2+\nu\xi_1,\xi_1\xi_2\cdot \nu(2+\nu)]$ is
equal to $(2+\nu)|\xi_1|\cdot |\xi_2(2+\nu)+\xi_1\nu|\approx
2^{k_1+k_2}$, see \eqref{jj11.1} and the definition of the set
$S''$.
\end{proof}

\newtheorem{Main10}[Main9]{Lemma}
\begin{Main10}\label{Main10}
Assume $k_1,k_2,k_3\in \Z$, $j_1,j_2,j_3\in\Z_+$, and
$f_i:\R^3\to\R_+$ are $L^2$ functions supported in
$D_{k_i,\infty,j_i}$, $i=1,2,3$. Then
\begin{equation}\label{mola1}
\int_{\R^3}(f_1\ast f_2)\cdot f_3\lesssim   2^{(j_1+j_2+j_3)/2}\cdot
2^{-\max(j_1,j_2,j_3)/2}\cdot
\|f_1\|_{L^2}\|f_2\|_{L^2}\|f_3\|_{L^2}.
\end{equation}
\end{Main10}

\begin{proof} Using the symmetry \eqref{mola6}, we may assume
$j_3=\max(j_1,j_2,j_3)$. Then
\begin{equation*}
\int_{\R^3}(f_1\ast f_2)\cdot f_3\lesssim   \|f_3\|_{L^2}\cdot
\|f_1\ast f_2\|_{L^2}\lesssim 
 \|f_3\|_{L^2}\|\mathcal{F}^{-1}(f_1)\|_{L^4}\|\mathcal{F}^{-1}(f_2)\|_{L^4}.
\end{equation*}
We use the scale-invariant Strichartz estimate of \cite{ArSa}:
\begin{equation}\label{io30}
\left\|\int_{\R^2}\phi(\xi,\mu)e^{ix\cdot\xi}e^{iy\cdot\mu}e^{it\cdot\omega(\xi,\mu)}\,d\xi
d\mu\right\|_{L^4_{x,y,t}}\lesssim   \|\phi\|_{L^2},
\end{equation}
for any $\phi\in L^2(\R^2)$. With $f_i^\#$, $i=1,2$, defined as in
the proof of Lemma \ref{Main9}, we estimate
\begin{equation*}
\begin{split}
&\left\|\int_{\R^3}f_i(\xi,\mu,\tau)\cdot e^{ix\cdot\xi}e^{iy\cdot\mu}e^{it\cdot\tau}\,d\xi d\mu d\tau\right\|_{L^4_{x,y,t}}\\
&=\left\|\int_{\R^3}f_i^\#(\xi,\mu,\theta)\cdot e^{it\cdot\theta}\cdot e^{ix\cdot\xi}e^{iy\cdot\mu}e^{it\cdot\omega(\xi,\mu)}\,d\xi d\mu d\theta\right\|_{L^4_{x,y,t}}\\
&\lesssim   2^{j_i/2}\|f_i^\#(\xi,\mu,\theta)\|_{L^2},
\end{split}
\end{equation*}
which gives \eqref{mola1}.
\end{proof}

As a consequence of Lemma \ref{Main9} and Lemma \ref{Main10}, we
have the following $L^2$ bilinear estimates.

\newtheorem{Main9co}[Main9]{Corollary}
\begin{Main9co}\label{Main9co}
(a) Assume $k_1,k_2,k\in \Z$, $j_1,j_2,j\in\Z_+$, and
$f_i:\R^3\to\R_+$ are $L^2$ functions supported in
$D_{k_i,\infty,j_i}$, $i=1,2$. Then
\begin{equation}\label{jj1}
\|\mathbf{1}_{D_{k,\infty,j}}\cdot(f_1\ast f_2)\|_{L^2}\lesssim 
 2^{(j_1+j_2+j)/2}
(2^{\max(j_1,j_2,j)}+2^{k_1+k_2+k})^{-1/2}\cdot
\|f_1\|_{L^2}\|f_2\|_{L^2}.
\end{equation}

(b) Assume $k_1,k_2,k\in \Z$,  $j_1,j_2,j\in\Z_+$,
and $f_i:\R^3\to\R_+$ are $L^2$ functions supported in
$D_{k_i,\infty,j_i}$, $i=1,2$. If $k_1\leq 100$ then
\begin{equation}\label{jj2lowk}
\|\mathbf{1}_{D_{k,\infty,j}}\cdot (f_1\ast f_2)\|_{L^2}\lesssim
 2^{[k_1+\min(k_1,k_2,k)+\min(j_1,j_2,j)]/2}\cdot
\| p(\xi_1,\mu_1)\cdot f_1\|_{L^2}\|f_2\|_{L^2}.
\end{equation}
If $k_1\geq -100$ then
\begin{equation}\label{jj2highk}
\|\mathbf{1}_{D_{k,\infty,j}}\cdot(f_1\ast f_2)\|_{L^2}\lesssim
 2^{[2k_1+\min(k_1,k_2,k)+\min(j_1,j_2,j)]/2}\cdot
\| p(\xi_1,\mu_1)\cdot f_1\|_{L^2}\|f_2\|_{L^2}.
\end{equation}
\end{Main9co}

\begin{proof}[Proof of Corollary \ref{Main9co}] Part (a) follows from \eqref{jj1.6} and \eqref{mola1}. For part (b), recall (see \eqref{eq-3}) that $p(\xi,\mu)=1+|\mu|/( |\xi|+|\xi|^2)$. To prove \eqref{jj2lowk} we decompose
\begin{equation*}
f_1=f_{1,k_1}+\sum_{l_1=k_1+1}f_{1,l_1}=f_1\cdot \eta_0(\mu_1/2^{k_1})+\sum_{l_1=k_1+1}^\infty f_1\cdot \chi_{l_1}(\mu_1).
\end{equation*}
Using \eqref{jj2.22}, the left-hand side of \eqref{jj2lowk} is dominated by
\begin{equation*}
\begin{split}
&\sum_{l_1=k_1}^\infty \|\mathbf{1}_{D_{k,\infty,j}}\cdot (f_{1,l_1}\ast f_2)\|_{L^2}\lesssim2^{[\min(k_1,k_2,k)+\min(j_1,j_2,j)]/2}\|f_2\|_{L^2}\sum_{l_1=k_1}^\infty 2^{l_1/2}\| f_{1,l_1}\|_{L^2}\\
&\lesssim 2^{[\min(k_1,k_2,k)+\min(j_1,j_2,j)]/2}\|f_2\|_{L^2}\cdot 2^{k_1/2}\| p(\xi_1,\mu_1)\cdot f_1\|_{L^2},
\end{split}
\end{equation*}
as desired. To prove \eqref{jj2highk} we decompose
\begin{equation*}
f_1=f_{1,2k_1}+\sum_{l_1=2k_1+1}f_{1,l_1}=f_1\cdot \eta_0(\mu_1/2^{2k_1})+\sum_{l_1=2k_1+1}^\infty f_1\cdot \chi_{l_1}(\mu_1).
\end{equation*}
Using \eqref{jj2.22}, the left-hand side of \eqref{jj2lowk} is dominated by
\begin{equation*}
\begin{split}
&\sum_{l_1=2k_1}^\infty \|\mathbf{1}_{D_{k,\infty,j}}\cdot (f_{1,l_1}\ast f_2)\|_{L^2}\lesssim2^{[\min(k_1,k_2,k)+\min(j_1,j_2,j)]/2}\|f_2\|_{L^2}\sum_{l_1=2k_1}^\infty 2^{l_1/2}\| f_{1,l_1}\|_{L^2}\\
&\lesssim 2^{[\min(k_1,k_2,k)+\min(j_1,j_2,j)]/2}\|f_2\|_{L^2}\cdot 2^{k_1}\| p(\xi_1,\mu_1)\cdot f_1\|_{L^2},
\end{split}
\end{equation*}
as desired.
\end{proof}

\section{Energy estimates}\label{energyproof}

In this section we prove the energy estimates in
Proposition~\ref{Lemmad1} and Proposition~\ref{Lemmak4}. To prove dyadic energy estimates we introduce a new Littlewood-Paley
decomposition with smooth symbols.
With
\[
\chi_k(\xi)=\eta_0(\xi/2^k)-\eta_0(\xi/2^{k-1}), \qquad k\in\Z,
\]
let $\widetilde{P}_k$ denote the operator on $L^2(\mathbb{R}^3)$
defined by the Fourier multiplier $(\xi,\mu,\tau)\to\chi_k(\xi)$.  By
a slight abuse of notation, we also let $\widetilde{P}_k$ denote the
operator on $L^2(\mathbb{R}^2)$ defined by the Fourier multiplier
$(\xi,\mu)\to\chi_k(\xi)$. For $l\in\Z$ let
\begin{equation*}
\widetilde{P}_{\leq l}=\sum_{k\leq l}\widetilde{P}_k,\qquad \widetilde{P}_{\geq l}=\sum_{k\geq l}\widetilde{P}_k.
\end{equation*}

Assume that, for some $k\in\Z$ and $u,v\in C([-T,T]:\overline{E}_k)$
\begin{equation}
\label{eq-linf}
\begin{cases}
\partial_tu+\partial_x^3u-\partial_x^{-1}\partial_y^2u=v\text{ on }\R^2\times (-T,T);\\
u(0)=\phi.
\end{cases}
\end{equation}
We multiply by $u$ and integrate to conclude that
\begin{equation} \label{supenk}
\sup_{|t_k|\leq T}\|u(t_k)\|_{L^2}^2 \leq  \|\phi\|_{L^2}^2 +\sup_{|t_k|\leq T} \left| \int_{\R^2\times[0,t_k]}u\cdot v\,dxdydt\right|.
\end{equation}

To prove Proposition~\ref{Lemmad1} and Proposition~\ref{Lemmak4} we
need to replace $v$ by the corresponding bilinear expressions. Thus we
need to estimate integrals of trilinear forms. However, instead of
direct estimates we seek to take advantage of the special form of the
nonlinearities. This allows us to place the derivative in the
nonlinearity on the lowest frequency factor. We summarize the main
dyadic estimates we need in Lemma \ref{tri} below.

\begin{lemma} \label{tri}
(a) Assume $T\in(0,1]$, $k_1,k_3,k_3 \in \Z$
with $\max\{k_1,k_2,k_3\} \geq 0$, and $u_i \in  \overline F_{k_i}(T)$, $i=1,2,3$. Assume in addition that $u_i \in  F_{k_i}(T)$ for some $i\in\{1,2,3\}$. Then
\begin{equation}\label{tria}
\left|  \int_{\R^2\times[0,T]}u_1u_2u_3\,dxdydt\right| \lesssim
 2^{-\min(k_1,k_2,k_3)/2}\prod_{i=1}^3\|u_{k_i}\|_{\overline F_{k_i}(T)}.
\end{equation}

(b) Assume $T\in (0,1]$, $k\in\Z_+$, $k_1\leq k-10$, $u\in \overline{\F}^0(T)$, and $v\in F_{k_1}(T)$. Then
\begin{equation}\label{triab}
\left|
\int_{\R^2\times[0,T]}\widetilde{P}_k(u)\widetilde{P}_k(\partial_xu\cdot
\widetilde{P}_{k_1}(v))\,dxdydt\right| \lesssim
2^{k_1/2}\|v\|_{\overline{F}_{k_1}(T)}\negmedspace\sum_{|k'-k|\leq
10}\|\widetilde{P}_{k'}(u)\|_{\overline{F}_{k'}(T)}^2.
\end{equation}
\end{lemma}

\begin{proof}[Proof of Lemma \ref{tri}] For part (a), we may assume that
$k_1 \leq k_2 \leq k_3$.
In order for the integral to be nontrivial we must also have
$|k_2 -k_3| \leq 4$. The integral in the left-hand side of \eqref{tria} converges absolutely, since one of the factors is in $F_{k}(T)$, thus bounded. We fix extensions $\widetilde{u}_i\in \overline{F}_{k_i}$ such that $\|\widetilde{u}_i\|_{\overline{F}_{k_i}}\leq 2\|u_i\|_{\overline{F}_{k_i}(T)}$, $i=1,2,3$. Let $\gamma:\mathbb{R}\to[0,1]$ denote a smooth function supported in $[-1,1]$ with the property that
\begin{equation*}
\sum_{n\in\Z}\gamma^3(x-n)\equiv 1,\qquad x\in\R.
\end{equation*}
The left-hand side of \eqref{tria} is dominated by
\begin{equation}\label{bin1}
\begin{split}
C\sum_{|n|\leq C2^{k_3}}&\Big|
  \int_{\R^2\times\R}(\gamma(2^{k_3}t-n)\mathbf{1}_{[0,T]}(t)\widetilde{u}_1) \\ &
 \times(\gamma(2^{k_3}t-n)\mathbf{1}_{[0,T]}(t)\widetilde{u}_2)\cdot (\gamma(2^{k_3}t-n)\mathbf{1}_{[0,T]}(t)\widetilde{u}_3)\,dxdydt\Big|.
\end{split}
\end{equation}

To estimate the integrals in \eqref{bin1} we  observe that, in view of \eqref{jj1}, if $k_1,k_2,k_3\in \Z$, $f_{k_i}\in
X_{k_i}$, $i=1,2,3$, and $|m|\leq 1$ then
\begin{equation}\label{sd31}
\begin{split}
&\Big|\int_{\R^3}\int_{\R^3}m(\xi,\xi_1)\cdot
f_3(-\xi,-\mu,-\tau)f_2(\xi-\xi_1,\mu-\mu_1,\tau-\tau_1)\\
&\times f_1(\xi_1,\mu_1,\tau_1)\,d\xi d\mu d\tau d\xi_1 d\mu_1
d\tau_1\Big|\lesssim (1+2^{k_1+k_2+k_3})^{-1/2}\Pi,
\end{split}
\end{equation}
where $\Pi=||f_1||_{X_{k_1}}||f_2||_{X_{k_2}}||f_3||_{X_{k_3}}$. In addition, as in \eqref{sp7.3}, if $I\subseteq\R$ is an interval, $k\in\Z$, $f_k\in X_k$, and $f_k^I=\mathcal{F}(\mathbf{1}_I(t)\cdot \mathcal{F}^{-1}(f_k))$ then
\begin{equation*}
\sup_{j\in\Z_+}2^{j/2}\|\eta_j(\tau-\omega(\xi,\mu))\cdot
f_k^I\|_{L^2}\lesssim \|f_k\|_{X_k}.
\end{equation*}
Thus, using \eqref{jj1} again, if $k_1,k_2,k_3\in \Z$, $f_{k_i}\in
X_{k_i}$, $i=1,2,3$, $I_i\subseteq\R$, $i=1,2,3$, are intervals, and $|m|\leq 1$ then
\begin{equation}\label{sd32}
\begin{split}
\Big|&\int_{\R^3}\int_{\R^3}m(\xi,\xi_1)\cdot
f_3^{I_3}(-\xi,-\mu,-\tau)f_2^{I_2}(\xi-\xi_1,\mu-\mu_1,\tau-\tau_1)\\
& f_1^{I_1}(\xi_1,\mu_1,\tau_1)\,d\xi d\mu d\tau d\xi_1 d\mu_1
d\tau_1\Big|\lesssim 
(1+2^{k_1+k_2+k_3})^{-1/2}\max(1,k_1,k_2,k_3)^3\Pi.
\end{split}
\end{equation}

We apply now the bound \eqref{sd32} up to $4$ times (for the
integers $n$ for which
$\gamma(2^{k_3}t-n)\mathbf{1}_{[0,T]}(t)\neq\gamma(2^{k_3}t-n)$)
and the bound \eqref{sd31} about $2^{k_3}$ times to bound the sum
in \eqref{bin1} by the right-hand side of \eqref{tria} (using also
\eqref{sp7.3}). This completes the proof of part (a).

For part (b), we observe first that the expression in the
left-hand side of \eqref{triab} is dominated by
\begin{equation}\label{bin2}
\begin{split}
&C\left|  \int_{\R^2\times[0,T]}\widetilde{P}_k(u)\cdot \widetilde{P}_k(\partial_xu)\cdot \widetilde{P}_{k_1}(v)\,dxdydt\right|\\
&+C\left|  \int_{\R^2\times[0,T]}\widetilde{P}_k(u)\cdot [\widetilde{P}_k(\partial_xu\cdot \widetilde{P}_{k_1}(v))-\widetilde{P}_k(\partial_xu)\cdot \widetilde{P}_{k_1}(v)]\,dxdydt\right|.
\end{split}
\end{equation}
We integrate by parts and use \eqref{tria} to conclude that
\begin{equation}\label{bin3}
\begin{split}
\left|  \int_{\R^2\times[0,T]}\widetilde{P}_k(u)\cdot
  \widetilde{P}_k(\partial_xu)\cdot
  \widetilde{P}_{k_1}(v)\,dxdydt\right|\lesssim 2^{k_1/2}\|\widetilde{P}_{k_1}(v)\|_{\overline{F}_{k_1}(T)}\cdot \|\widetilde{P}_k(u)\|_{\overline{F}_k(T)}^2,
\end{split}
\end{equation}
which suffices for \eqref{triab}.

To control the term in the second line of \eqref{bin2} we fix extensions $\widetilde{u}$ of $u$ and $\widetilde{v}$ of $v$ and use the formula
\begin{equation}\label{sf1}
\begin{split}
&\mathcal{F}[\widetilde{P}_k(\widetilde{P}_{k_1}(\widetilde{v})\cdot \partial_x\widetilde{u})-\widetilde{P}_{k_1}(\widetilde{v})\cdot \widetilde{P}_k(\partial_x\widetilde{u})](\xi,\mu,\tau)\\
&=C\int_{\R^3}\mathcal{F}(\widetilde{P}_{k_1}(\partial_x\widetilde{v}))(\xi_1,\mu_1,\tau_1)\cdot \mathcal{F}(\widetilde{u})(\xi-\xi_1,\mu-\mu_1,\tau-\tau_1)\cdot m(\xi,\xi_1)\,d\xi_1d\mu_1 d\tau_1,
\end{split}
\end{equation}
where
\begin{equation*}
|m(\xi,\xi_1)|=\Big|\frac{(\xi-\xi_1)(\chi_k(\xi)-\chi_k(\xi-\xi_1))}{\xi_1}\Big|\lesssim
\sum_{|k'-k|\leq 4}\chi_{k'}(\xi-\xi_1).
\end{equation*}
The bound \eqref{triab} follows by decomposing the integral in the
second line of \eqref{bin2} into at most $C2^k$ integrals over
time-intervals of length $\approx 2^{-k}$ (as in \eqref{bin1}),
and using the formula \eqref{sf1} and the bounds \eqref{sd31} and
\eqref{sd32} to bound these integrals.
\end{proof}

We prove now Proposition \ref{Lemmad1} and Proposition \ref{Lemmak4}.

\begin{proof}[Proof of Proposition \ref{Lemmad1}]
Recall that $u$ solves the initial-value problem
\begin{equation}\label{bin10}
\begin{cases}
\partial_tu+\partial_x^3u-\partial_x^{-1}\partial_y^2u+\partial_x(u^2/2)=0\text{ on }\R^2\times[-T,T];\\
u(0)=\phi.
\end{cases}
\end{equation}
We observe that
\begin{equation}\label{sd2}
\begin{split}
&\|u\|^2_{\B^\sigma(T)}-\|P_{\leq 0}(\phi)\|_{\E^\sigma}^2\\
&\lesssim \sum_{k\geq 0}\sup_{t_k\in[-T,T]}(2^{2\sigma
k}\|\widetilde{P}_k(u(t_k))\|_{L^2}^2
+2^{(2\sigma-2)k}||\widetilde{P}_k(\partial_x^{-1}\partial_yu(t_k))||^2_{L^2}).
\end{split}
\end{equation}
Therefore it suffices to prove that for $\sigma\in\{1,2,3\}$
\begin{equation}\label{sd9}
\begin{split}
\sum_{k\geq 0}&\sup_{t_k\in[-T,T]}2^{2\sigma
k}\|\widetilde{P}_k(u(t_k))\|_{L^2}^2 +\sum_{k\geq
0}\sup_{t_k\in[-T,T]}2^{(2\sigma-2)k}
||\widetilde{P}_k(\partial_x^{-1}\partial_yu(t_k))||^2_{L^2})\\
&\lesssim \|\phi\|_{\E^\sigma}^2+ \|u\|_{\F^1(T)}\cdot
\|u\|_{\F^\sigma(T)}^2.
\end{split}
\end{equation}

We show first that
\begin{equation}\label{sd5}
\sum_{k\geq 0}\sup_{t_k\in[-T,T]}2^{2\sigma
k}\|\widetilde{P}_k(u(t_k))\|_{L^2}^2 -2^{2\sigma
k}\|\widetilde{P}_k(\phi)\|_{L^2}^2\lesssim \|u\|_{\F^1(T)}\cdot
\|u\|_{\F^\sigma(T)}^2.
\end{equation}
For $k\in\Z_+$ we use \eqref{supenk} and the equation
\eqref{bin10} to estimate the increment
\begin{equation}\label{sd10}
2^{2\sigma k}\|\widetilde{P}_k(u(t_k))\|_{L^2}^2-2^{2\sigma
k}\|\widetilde{P}_k(\phi)\|_{L^2}^2\lesssim 2^{2\sigma
 k}\Big|\int_{\R^2\times[0,{t_k}]}\widetilde{P}_k(u)\widetilde{P}_k(u\cdot
\partial_xu)\,dxdydt\Big|.
\end{equation}
The right-hand side of \eqref{sd10} is
dominated by
\begin{equation}\label{sd30}
\begin{split}
&C2^{2\sigma k}\sum_{k_1\leq
k-10}\Big|\int_{\R^2\times[0,{t_k}]}\widetilde{P}_k(u)\cdot
\widetilde{P}_k(\widetilde{P}_{k_1}(u)\cdot
\partial_xu)\,dxdydt\Big|\\
&+C2^{2\sigma k}\sum_{k_1\geq k-9,k_2\in\Z}\Big|\int_{\R^2\times[0,{t_k}]}\widetilde{P}_k^2(u)\cdot
\widetilde{P}_{k_1}(u)\cdot
\partial_x\widetilde{P}_{k_2}(u)\,dxdydt\Big|.
\end{split}
\end{equation}
Using \eqref{triab}, the sum in the first line of \eqref{sd30} is
dominated by
\begin{equation*}
\begin{split}
C\|{u}\|_{\F^1(T)}\cdot \sum_{|k'-k|\leq 10}2^{2\sigma
k'}\|\widetilde{P}_{k'}(u)\|_{\overline{F}_{k'}(T)}^2.
\end{split}
\end{equation*}
Using \eqref{tria}, the sum in the second line of \eqref{sd30} is
dominated by
\begin{equation*}
\begin{split}
& C2^{2\sigma k}\sum_{|k_1-k|\leq 10,k_2\leq k+10}
2^{k_2/2}\|\widetilde{P}_k(u)\|_{\overline{F}_k(T)}
\|\widetilde{P}_{k_1}(u)\|_{\overline{F}_{k_1}(T)}
\|\widetilde{P}_{k_2}(u)\|_{\overline{F}_{k_2}(T)}\\
&+C2^{2\sigma k}\sum_{k_1\geq k+10, |k_2-k_1|\leq 10}
2^{k_2-k/2}\|\widetilde{P}_k(u)\|_{\overline{F}_k(T)}
\|\widetilde{P}_{k_1}(u)\|_{\overline{F}_{k_1}(T)}
\|\widetilde{P}_{k_2}(u)\|_{\overline{F}_{k_2}(T)}\\
&\lesssim \|{u}\|_{\F^1(T)}\cdot \sum_{|k'-k|\leq 20}2^{2\sigma
k'}\|\widetilde{P}_{k'}(u)\|_{\overline{F}_{k'}(T)}^2+2^{k/2}
\|\widetilde{P}_k(u)\|_{\overline{F}_k(T)}\cdot\|u\|_{\F^\sigma(T)}^2.
\end{split}
\end{equation*}
The bound \eqref{sd5} follows.

We show now that
\begin{equation}\label{sq5}
\begin{split}
\sum_{k\geq 0}\sup_{t_k\in[-T,T]}2^{(2\sigma-2)k}
\|\widetilde{P}_k(\partial_x^{-1}\partial_yu(t_k))\|_{L^2}^2 -2^{(2\sigma-2) k}\|\widetilde{P}_k(\partial_x^{-1}\partial_y\phi)\|_{L^2}^2\\
\lesssim \|u\|_{\F^1(T)}\cdot \|u\|_{\F^\sigma(T)}^2.
\end{split}
\end{equation}
For $k\in\Z_+$ and $t_k\in[-T,T]$ we use \eqref{supenk} and the
the equation \eqref{bin10} to estimate the increment
\begin{equation}\label{sq10}
\begin{split}
2^{(2\sigma-2)
k}\|\widetilde{P}_k(\partial_x^{-1}\partial_yu(t_k))\|_{L^2}^2-2^{(2\sigma
-2)k}\|\widetilde{P}_k(\partial_x^{-1}\partial_y\phi)\|_{L^2}^2\\
\lesssim 2^{(2\sigma-2)
 k}\Big|\int_{\R^2\times[0,{t_k}]}\widetilde{P}_k(\partial_x^{-1}\partial_yu)\widetilde{P}_k(u\cdot
\partial_yu)\,dxdydt\Big|.
\end{split}
\end{equation}
The right-hand side of \eqref{sq10} is dominated by
\begin{equation}\label{sq30}
\begin{split}
&C2^{(2\sigma-2) k}\sum_{k_1\leq
k-10}\Big|\int_{\R^2\times[0,{t_k}]}\widetilde{P}_k(v)\cdot
\widetilde{P}_k(\widetilde{P}_{k_1}(u)\cdot
\partial_xv)\,dxdydt\Big|\\
&+C2^{(2\sigma-2) k}\sum_{k_1\geq
k-9,k_2\in\Z}\Big|\int_{\R^2\times[0,{t_k}]}\widetilde{P}_k^2(v)\cdot
\widetilde{P}_{k_1}(u)\cdot
\partial_x\widetilde{P}_{k_2}(v)\,dxdydt\Big|,
\end{split}
\end{equation}
where $v=\partial_x^{-1}\partial_yu$. Using \eqref{triab}, the sum
in the first line of \eqref{sq30} is dominated by
\begin{equation*}
\begin{split}
C\|{u}\|_{\F^1(T)}\cdot \sum_{|k'-k|\leq 10}2^{(2\sigma-2)
k'}\|\widetilde{P}_{k'}(\partial_x^{-1}\partial_yu)\|_{\overline{F}_{k'}(T)}^2.
\end{split}
\end{equation*}
Using \eqref{tria}, the sum in the second line of \eqref{sq30} is
dominated by
\begin{equation*}
\begin{split}
&C2^{(2\sigma-2) k}\sum_{|k_1-k|\leq 10,k_2\leq k+10}
2^{k_2/2}\|\widetilde{P}_k(v)\|_{\overline{F}_k(T)}
\|\widetilde{P}_{k_1}(u)\|_{\overline{F}_{k_1}(T)}
\|\widetilde{P}_{k_2}(v)\|_{\overline{F}_{k_2}(T)}\\
&+C2^{(2\sigma-2) k}\sum_{k_1\geq k+10, |k_2-k_1|\leq 10}
2^{k_2-k/2}\|\widetilde{P}_k(v)\|_{\overline{F}_k(T)}
\|\widetilde{P}_{k_1}(u)\|_{\overline{F}_{k_1}(T)}
\|\widetilde{P}_{k_2}(v)\|_{\overline{F}_{k_2}(T)}\\
&\lesssim \|{u}\|_{\F^1(T)}\cdot \sum_{|k'-k|\leq 20}2^{2\sigma
k'}\|\widetilde{P}_{k'}(u)\|_{{F}_{k'}(T)}^2+C2^{k/2}
\|\widetilde{P}_k(u)\|_{{F}_k(T)}\cdot\|u\|_{\F^\sigma(T)}^2.
\end{split}
\end{equation*}
The bound \eqref{sq5} follows, which completes the proof of
Proposition \ref{Lemmad1}.
\end{proof}

\begin{proof}[Proof of Proposition \ref{Lemmak4}] Recall that
$u=P_{\geq -10}(u)$ solves the equation
\begin{equation}\label{bin20}
\begin{cases}
\partial_tu+\partial_x^3u-\partial_x^{-1}\partial_y^2u=P_{\geq -10}(v\cdot \partial_xu)+\sum\limits_{m=1}^3P_{\geq -10}(w_m\cdot w'_m)+P_{\geq -10}(h);\\
u(0)=\phi,
\end{cases}
\end{equation}
on $\R^2\times(-T,T)$. It suffices to prove that
\begin{equation}\label{su10}
\begin{split}
&\sum_{k\geq 0}\sup_{t_k\in[-T,T]}\big( \|\widetilde{P}_k(u(t_k))\|_{L^2}^2-\|\widetilde{P}_k(\phi)\|_{L^2}^2\big) \\
&\lesssim \|v\|_{\F^1(T)}\cdot
\|u\|_{\overline{\F}^0(T)}^2+ \sum_{m=1}^3\|u\|_{\overline{\F}^0(T)}
\|w_m\|_{\overline{\F}^0(T)}\|w'_m\|_{\overline{\F}^0(T)}.
\end{split}
\end{equation}
Using \eqref{supenk} and the equation \eqref{bin20}, for $k\geq 0$
\begin{equation}\label{su11}
\begin{split}\!\!\! \|\widetilde{P}_k(u(t_k))\|_{L^2}^2-&\|\widetilde{P}_k(\phi)\|_{L^2}^2
\lesssim\Big|\int_{\R^2\times[0,t_k]}\widetilde{P}_k(u)\cdot
\widetilde{P}_k(\widetilde{P}_{\leq k-10}(v)\cdot \partial_xu)\,dxdydt\Big|\\
&+\Big|\int_{\R^2\times[0,t_k]}\widetilde{P}_k^2(u)\cdot
\partial_xu\cdot \widetilde{P}_{\geq {k-9}}(v)\,dxdydt\Big|\\
&+\sum_{m=1}^3\Big|\int_{\R^2\times[0,t_k]}\widetilde{P}_k^2(u)\cdot
w_m\cdot w'_m\,dxdydt\Big|.
\end{split}
\end{equation}
We observe that the term $P_{\geq -10}(h)$ plays no role in the
proof of \eqref{su10} (this term is needed, however, to prove the bounds \eqref{amy50} and
\eqref{amy51}).

Using \eqref{triab},
\begin{equation*}
\sum_{k\geq 0}\Big|\int_{\R^2\times[0,t_k]}\widetilde{P}_k(u)\cdot
\widetilde{P}_k(\widetilde{P}_{\leq k-10}(v)\cdot
\partial_xu)\,dxdydt\Big|\lesssim \|v\|_{\F^1(T)}\|u\|_{\overline{\F}^0}^2.
\end{equation*}
Using \eqref{tria},
\begin{equation*}
\begin{split}
&\sum_{k\geq
0}\Big|\int_{\R^2\times[0,t_k]}\widetilde{P}_k^2(u)\cdot
\partial_xu\cdot \widetilde{P}_{\geq {k-9}}(v)\,dxdydt\Big|\\
&\lesssim \sum_{k\geq 0}\sum_{k_2\geq k-9}\sum_{k_1\leq k_2+20}
\Big|\int_{\R^2\times[0,t_k]}\widetilde{P}_k^2(u)\cdot
\partial_x\widetilde{P}_{k_1}(u)\cdot \widetilde{P}_{k_2}(v)\,dxdydt\Big|\\
&\lesssim \sum_{k\geq 0}\sum_{k_2\geq k-9}\sum_{k_1\leq k_2+20}
2^{k_1-\min(k_1,k)/2}\cdot
\|\widetilde{P}_k(u)\|_{\overline{F}_k(T)}
\|\widetilde{P}_{k_1}(u)\|_{\overline{F}_{k_1}(T)}
\|\widetilde{P}_{k_2}(v)\|_{\overline{F}_{k_2}(T)}\\
&\lesssim \|v\|_{\F^1(T)}\|u\|_{\overline{\F}^0}^2.
\end{split}
\end{equation*}
Using \eqref{tria}, with
$k_\mathrm{med}=k+k_1+k_2-\min(k,k_1,k_2)-\max(k,k_1,k_2)$,
$k_{\max}=\max(k,k_1,k_2)$
\begin{equation*}
\begin{split}
&\sum_{k\geq
0}\Big|\int_{\R^2\times[0,t_k]}\widetilde{P}_k^2(u)\cdot w_m\cdot
w'_m\,dxdydt\Big|\\
&\lesssim
\sum_{k,k_1,k_2\in\Z}\Big|\int_{\R^2\times[0,t_k]}\widetilde{P}_k^2(u)\cdot
\widetilde{P}_{k_1}(w_m)\cdot
\widetilde{P}_{k_2}(w'_m)\,dxdydt\Big|\\
&\lesssim \sum_{|k_\mathrm{med}-k_{\max}|\leq
10}2^{-\min(k,k_1,k_2)/2}\|\widetilde{P}_k(u)\|_{\overline{F}_k(T)}
\|\widetilde{P}_{k_1}(w_m)\|_{\overline{F}_{k_1}(T)}
\|\widetilde{P}_{k_2}(w'_m)\|_{\overline{F}_{k_2}(T)}\\
&\lesssim \|u\|_{\overline{\F}^0(T)}
\|w_m\|_{\overline{\F}^0(T)}\|w'_m\|_{\overline{\F}^0(T)}.
\end{split}
\end{equation*}
This last inequality uses the fact that, for any
$v\in\overline{\F}^0(T)$,
\[
\sum_{k\in\Z}2^{-k/2}\|\widetilde{P}_k(v)\|_{\overline{F}_k(T)}\lesssim
\|v\|_{\overline{\F}^0(T)},
\]
 which is the main reason for the
low-frequency condition on functions in $\overline{\F}^0$.

The main bound \eqref{su10} follows, which completes the proof of
the proposition.
\end{proof}

\section{Dyadic bilinear estimates, I}\label{dyadic1}

In this section we prove several dyadic bounds which are used in
the proof of Proposition \ref{Lemmab1} (a). We estimate first $\mathrm{Low}\times\mathrm{High}\to\mathrm{High}$ interactions.

\newtheorem{Lemmab2}{Lemma}[section]
\begin{Lemmab2}\label{Lemmab2}
Assume $k,k_1,k_2\in\Z$, $k_1\leq k_2$, $k_1\leq 0$, $k\geq 0$,
$|k_2-k|\leq 40$, $u_{k_1}\in F_{k_1}$ and $v_{k_2}\in F_{k_2}$.
Then
\begin{equation}\label{on10}
\|P_k(\partial_x(u_{k_1}v_{k_2}))\|_{N_k}\lesssim 2^{k_1}\cdot
\|u_{k_1}\|_{F_{k_1}}\cdot \|v_{k_2}\|_{F_{k_2}}.
\end{equation}
\end{Lemmab2}

\begin{proof}[Proof of Lemma \ref{Lemmab2}] Using the definitions and \eqref{sb21}, the left-hand side of \eqref{on10} is dominated by
\begin{equation*}
\begin{split}
C\sup_{t_k\in\R}\|p(\xi,\mu)&(\tau-\omega(\xi,\mu)+i2^k)^{-1}\cdot 2^k\mathbf{1}_{I_k}(\xi)\cdot\\
&\mathcal{F}[u_{k_1}\cdot
\eta_0(2^{k}(t-t_k))]\ast\mathcal{F}[v_{k_2}\cdot
\eta_0(2^{k}(t-t_k))]\|_{X_k}.
\end{split}
\end{equation*}
Let $f_{k_1}=\mathcal{F}[u_{k_1}\cdot \eta_0(2^{k}(t-t_k))]$ and
$f_{k_2}=\mathcal{F}[v_{k_2}\cdot \eta_0(2^{k}(t-t_k))]$. Using
the bounds \eqref{sp7.2} and \eqref{sb21}, it suffices to prove
that if $j_1,j_2\geq k$, and
$f_{k_i,j_i}:\mathbb{R}^3\to\mathbb{R}_+$ are supported in
$D_{k_i,\infty,j_i}$, $i=1,2$, then
\begin{equation}\label{on11}
\begin{split}
2^k&\sum_{j\geq k}2^{-j/2}\|\mathbf{1}_{D_{k,\infty,j}}\cdot
p(\xi,\mu)
\cdot (f_{k_1,j_1}\ast f_{k_2,j_2})\|_{L^2}\\
&\lesssim 2^{k_1}\cdot 2^{j_1/2}\|p(\xi_1,\mu_1)\cdot
f_{k_1,j_1}\|_{L^2}\cdot 2^{j_2/2}\|p(\xi_2,\mu_2)\cdot
f_{k_2,j_2}\|_{L^2}.
\end{split}
\end{equation}
Since $j,j_1,j_2 \geq k$ it suffices to prove the $L^2$ product estimate
\begin{equation}\label{on11a}
\begin{split}
 &  \| \mathbf{1}_{D_{k,\infty,j}}\cdot p(\xi,\mu)
  \cdot (f_{k_1,j_1}\ast f_{k_2,j_2})\|_{L^2}\\
  &\lesssim 2^{k_1}\cdot  2^{\min(j_1,j_2)/2}\|p(\xi_1,\mu_1)\cdot
  f_{k_1,j_1}\|_{L^2}\cdot \|p(\xi_2,\mu_2)\cdot
  f_{k_2,j_2}\|_{L^2}.
\end{split}
\end{equation}
Using the obvious bound
\begin{equation} \label{pxieta}
p(\xi,\mu) \lesssim |\xi_1| |\xi |^{-2} p(\xi_1,\mu_1) + p(\xi_2,\mu_2),
\end{equation}
this is a consequence  the estimates
\begin{equation*}
\begin{split}
 &  \|  \mathbf{1}_{D_{k,\infty,j}}\cdot (f_{k_1,j_1}\ast f_{k_2,j_2})\|_{L^2}
  \lesssim 2^{k_1/2 + k_2} 2^{\min(j_1,j_2)/2}
\|  f_{k_1,j_1}\|_{L^2} \|p(\xi_2,\mu_2)
  f_{k_2,j_2}\|_{L^2},
\end{split}
\end{equation*}
and
\begin{equation*}
\begin{split}
  &  \| \mathbf{1}_{D_{k,\infty,j}}\cdot (f_{k_1,j_1}\ast f_{k_2,j_2})\|_{L^2}
 \lesssim 2^{k_1}\cdot  2^{\min(j_1,j_2)/2}\|p(\xi_1,\mu_1)
  f_{k_1,j_1}\|_{L^2}\cdot \|
  f_{k_2,j_2}\|_{L^2},
\end{split}
\end{equation*}
which follow from \eqref{jj2highk} and \eqref{jj2lowk} repectively.
\end{proof}

\newtheorem{Lemmab3}[Lemmab2]{Lemma}
\begin{Lemmab3}\label{Lemmab3}
Assume $k,k_1,k_2\in\Z$, $k_1\leq k_2$, $k_1\geq 0$, $k\geq 0$,
$|k_2-k|\leq 40$, $u_{k_1}\in F_{k_1}$, and $v_{k_2}\in F_{k_2}$.
Then
\begin{equation}\label{on30}
\|P_k(\partial_x(u_{k_1}v_{k_2}))\|_{N_k}\lesssim (1+k_1) 2^{-k_1/2}\cdot
\|u_{k_1}\|_{F_{k_1}}\cdot \|v_{k_2}\|_{F_{k_2}}.
\end{equation}
\end{Lemmab3}

\begin{proof}[Proof of Lemma \ref{Lemmab3}] As in
  the proof of Lemma \ref{Lemmab2}, using the definitions and the
  bounds \eqref{sp7.2} and \eqref{sb21}, it suffices to prove that if
  $j_1,j_2\geq k$ and $f_{k_i,j_i}:\mathbb{R}^3\to\mathbb{R}_+$ are supported in
  $D_{k_i,\infty,j_i}$, $i=1,2$, then
\begin{equation}\label{on21}
\begin{split}
2^k&\sum_{j\geq k}2^{-j/2}\|\mathbf{1}_{D_{k,\infty,j}}\cdot p(\xi,\mu)\cdot (f_{k_1,j_1}\ast f_{k_2,j_2})\|_{L^2}\\
&\lesssim (1+k_1) 2^{-k_1/2}\cdot 2^{j_1/2}\|p(\xi_1,\mu_1)\cdot
f_{k_1,j_1}\|_{L^2}\cdot 2^{j_2/2}\|p(\xi_2,\mu_2)\cdot
f_{k_2,j_2}\|_{L^2}.
\end{split}
\end{equation}
Since $j,j_1,j_2 \geq k$, the large modulations $j \geq  k+4k_1$ in the output
are controlled by the $L^2$ product estimate
\begin{equation}\label{on11a2}
\begin{split}
 &  \| \mathbf{1}_{D_{k,\infty,j}}\cdot p(\xi,\mu)
  \cdot (f_{k_1,j_1}\ast f_{k_2,j_2})\|_{L^2}\\
  &\lesssim 2^{3k_1/2}\cdot  2^{\min(j_1,j_2)/2} \|p(\xi_1,\mu_1)\cdot
  f_{k_1,j_1}\|_{L^2}\cdot \|p(\xi_2,\mu_2)\cdot
  f_{k_2,j_2}\|_{L^2}.
\end{split}
\end{equation}
In this case we have
\begin{equation} \label{pxieta2}
p(\xi,\mu) \lesssim |\xi_1|^2 |\xi |^{-2} p(\xi_1,\mu_1) + p(\xi_2,\mu_2).
\end{equation}
Hence \eqref{on11a2} is a consequence of the estimates
\begin{equation*}
\begin{split}
 &  \| \mathbf{1}_{D_{k,\infty,j}}\cdot(f_{k_1,j_1}\ast f_{k_2,j_2})\|_{L^2}
  \lesssim 2^{k_1/2+k_2} 2^{\min(j_1,j_2)/2}
\|  f_{k_1,j_1}\|_{L^2} \|p(\xi_2,\mu_2)
  f_{k_2,j_2}\|_{L^2}
\end{split}
\end{equation*}
and
\begin{equation*}
\begin{split}
  &  \| \mathbf{1}_{D_{k,\infty,j}}\cdot(f_{k_1,j_1}\ast f_{k_2,j_2})\|_{L^2}
 \lesssim 2^{3k_1/2} 2^{\min(j_1,j_2)/2}\|p(\xi_1,\mu_1)
  f_{k_1,j_1}\|_{L^2}\cdot \|
  f_{k_2,j_2}\|_{L^2},
\end{split}
\end{equation*}
 which follow from \eqref{jj2highk}.

 It remains to estimate the small modulations $k \leq j \leq k+4k_1$ in
 the output.  There are about $1+k_1$ possible values for $j$,
 therefore we need to prove that
\begin{equation}\label{on21af}
\begin{split}
&\|\mathbf{1}_{D_{k,\infty,j}}\cdot p(\xi,\mu)\cdot (f_{k_1,j_1}\ast f_{k_2,j_2})\|_{L^2}\\
&\lesssim 2^{-k_1/2-k}\cdot 2^{(j+j_1+j_2)/2}\|p(\xi_1,\mu_1)\cdot
f_{k_1,j_1}\|_{L^2}\cdot \|p(\xi_2,\mu_2)\cdot
f_{k_2,j_2}\|_{L^2}.
\end{split}
\end{equation}
We observe that
\[
\left( \frac{\mu_2}{\xi_2} - \frac{\mu}{\xi} \right)^2 \lesssim
 |\xi_1|^2 + |\xi_1| |\xi|^{-2} |\Omega((\xi_1,\mu_1),(\xi_2,\mu_2))|
\]
which leads to
\begin{equation} \label{etapexi0}
 \frac{|\mu|}{|\xi|} \lesssim  \frac{|\mu_2|}{|\xi_2|} + |\xi_2| +
  |\xi_1|^\frac12 |\xi|^{-1} |\Omega((\xi_1,\mu_1),(\xi_2,\mu_2))|^{1/2},
\end{equation}
therefore
\begin{equation} \label{etaxibc}
p(\xi,\mu) \lesssim  p(\xi_2,\mu_2) + 2^{k_1/2} 2^{-2k_2}
2^{\max(j_1,j_2,j)/2}.
\end{equation}
We eliminate the expression  $p(\xi,\mu)$  on the left using
\eqref{etaxibc}, neglecting the remaining $ p(\xi,\mu)$
factors on the right. Then it suffices to show that
\[
\begin{split}
&\|\mathbf{1}_{D_{k,\infty,j}} \cdot (f_{k_1,j_1}\ast f_{k_2,j_2})\|_{L^2}
\lesssim 2^{-k_1/2-k} 2^{(j+j_1+j_2)/2}\|
f_{k_1,j_1}\|_{L^2} \|
f_{k_2,j_2}\|_{L^2},
\end{split}
\]
and
\[
\begin{split}
&\|\mathbf{1}_{D_{k,\infty,j}}\cdot (f_{k_1,j_1}\ast
f_{k_2,j_2})\|_{L^2} \lesssim 2^{-k_1+k_2} 2^{(j+j_1+j_2)/2} 2^{-\max(j_1,j_2,j)/2}
\|
f_{k_1,j_1}\|_{L^2} \|
f_{k_2,j_2}\|_{L^2}.
\end{split}
\]
Both these estimates follow from \eqref{jj1}.
\end{proof}

We estimate now $\mathrm{High}\times\mathrm{High}\to\mathrm{Low}$ interactions. Let
$\gamma:\mathbb{R}\to[0,1]$ denote a smooth function supported in
$[-1,1]$ with the property that $\sum_{m\in\Z}\gamma^2(x-m)\equiv
1$, $x\in\R$.

\newtheorem{Lemmab4}[Lemmab2]{Lemma}
\begin{Lemmab4}\label{Lemmab4}
Assume $k,k_1,k_2\in\Z_+$, $|k_1-k_2|\leq 4$, $k\leq
\min(k_1,k_2)-30$, $u_{k_1}\in F_{k_1}$, and $v_{k_2}\in F_{k_2}$.
Then
\begin{equation}\label{on80}
\|P_k(\partial_x(u_{k_1}v_{k_2}))\|_{N_k}\lesssim k_2
 2^{k_2-3k/2} \cdot
\|u_{k_1}\|_{F_{k_1}}\cdot \|v_{k_2}\|_{F_{k_2}}.
\end{equation}
\end{Lemmab4}

\begin{proof}[Proof of Lemma \ref{Lemmab4}] Using the definitions and \eqref{sb21}, the left-hand side of \eqref{on80} is dominated by
\begin{equation*}
\begin{split}
&C\sup_{t_k\in\R}\Big\|p(\xi,\mu)(\tau-\omega(\xi,\mu)+i2^k)^{-1}\cdot 2^k\mathbf{1}_{I_k}(\xi)\cdot\sum_{|m|\leq C2^{k_2-k}}\\
&\mathcal{F}[u_{k_1}\eta_0(2^{k}(t-t_k))\gamma(2^{k_2}(t-t_k)-m)]\negmedspace\ast\negmedspace\mathcal{F}[v_{k_2}\eta_0(2^{k}(t-t_k))\gamma(2^{k_2}(t-t_k)-m)]\Big\|_{X_k}
\end{split}
\end{equation*}
Using the definitions and the bounds \eqref{sp7.2} and
\eqref{sb21}, it suffices to prove that if $j_1, j_2\geq k_2$, and
$f_{k_i,j_i}:\mathbb{R}^3\to\mathbb{R}_+$ are supported in
$D_{k_i,\infty,j_i}$, $i=1,2$, then
\begin{equation}\label{on81}
\begin{split}
&2^k2^{k_2-k}\sum_{j\geq k}2^{-j/2}\|\mathbf{1}_{D_{k,\infty,j}}\cdot p(\xi,\mu)\cdot (f_{k_1,j_1}\ast f_{k_2,j_2})\|_{L^2}\\
&\lesssim k_2 2^{k_2-3k/2}\cdot 2^{j_1/2}\|p(\xi_1,\mu_1)\cdot
f_{k_1,j_1}\|_{L^2}\cdot 2^{j_2/2}\|p(\xi_2,\mu_2)\cdot
f_{k_2,j_2}\|_{L^2}.
\end{split}
\end{equation}
Due to the rough estimate
\begin{equation} \label{pximu1}
p(\xi,\mu) \lesssim 2^{2k_2 -2k} (p(\xi_1,\mu_1)+p(\xi_2,\mu_2))
\end{equation}
the bound \eqref{on81} follows from \eqref{jj1} in the region for $j/2 \geq 2k_2 - k/2$. Therefore
it remains to prove that
\begin{equation}\label{on81a}
\begin{split}
&2^{-j/2}\|\mathbf{1}_{D_{k,\infty,j}}\cdot p(\xi,\mu)\cdot (f_{k_1,j_1}\ast f_{k_2,j_2})\|_{L^2}\\
&\lesssim 2^{-3k/2}\cdot 2^{j_1/2}\|p(\xi_1,\mu_1)\cdot
f_{k_1,j_1}\|_{L^2}\cdot 2^{j_2/2}\|p(\xi_2,\mu_2)\cdot
f_{k_2,j_2}\|_{L^2}.
\end{split}
\end{equation}
We now seek to improve \eqref{pximu1}. We observe that
\[
\left( \frac{\mu_1}{\xi_1} - \frac{\mu}{\xi} \right)^2 \lesssim
 |\xi_2|^2 + |\xi|^{-1} |\Omega((\xi_1,\mu_1),(\xi_2,\mu_2))|,
\]
which leads to
\begin{equation} \label{etapexi}
 \frac{|\mu|}{|\xi|} \lesssim  \frac{|\mu_1|}{|\xi_1|} + |\xi_2| +
  |\xi|^{-1/2} |\Omega((\xi_1,\mu_1),(\xi_2,\mu_2))|,
\end{equation}
therefore
\[
p(\xi,\mu) \lesssim 2^{k_2-k}  p(\xi_1,\mu_1) + 2^{-3k/2}
2^{\max(j_1,j_2,j)/2}
\]
Thus \eqref{on81a} follows from the bounds
\begin{equation}\label{on81b}
\begin{split}
&\|\mathbf{1}_{D_{k,\infty,j}} \cdot (f_{k_1,j_1}
\ast f_{k_2,j_2})\|_{L^2} \lesssim 2^{-k_2-k/2} 2^{(j+j_1+j_2)/2}\|
f_{k_1,j_1}\|_{L^2}\|
f_{k_2,j_2}\|_{L^2},
\end{split}
\end{equation}
and
\begin{equation}\label{on81c}
\begin{split}
&\|\mathbf{1}_{D_{k,\infty,j}}\cdot  (f_{k_1,j_1}
\ast f_{k_2,j_2})\|_{L^2} \lesssim
2^{(j+j_1+j_2)/2 -\max(j,j_1,j_2)/2}
\|f_{k_1,j_1}\|_{L^2}\|
f_{k_2,j_2}\|_{L^2},
\end{split}
\end{equation}
both of which are consequences of \eqref{jj1}.
\end{proof}

\newtheorem{Lemmab5}[Lemmab2]{Lemma}
\begin{Lemmab5}\label{Lemmab5}
Assume that $k_1,k_2\in\Z_+$, $|k_1-k_2|\leq 4$,
$k\in\Z\cap(-\infty,0]$, $k\leq \min(k_1,k_2)-30$, $u_{k_1}\in
F_{k_1}$, and $v_{k_2}\in F_{k_2}$. Then
\begin{equation}\label{on280}
\|P_k(\partial_x(u_{k_1}v_{k_2}))\|_{N_k}\lesssim (k_2-k) 2^{k_2+k/2}
\cdot \|u_{k_1}\|_{F_{k_1}}\cdot \|v_{k_2}\|_{F_{k_2}}.
\end{equation}
\end{Lemmab5}

\begin{proof}[Proof of Lemma \ref{Lemmab5}] As in the proof of
Lemma \ref{Lemmab4}, using the definitions and the bounds
\eqref{sp7.2} and \eqref{sb21}, it suffices to prove that if $j_1,
j_2\geq k_2$, and $f_{k_i,j_i}:\mathbb{R}^3\to\mathbb{R}_+$
are supported in $D_{k_i,\infty,j_i}$, $i=1,2$, then
\begin{equation}\label{on281}
\begin{split}
&2^k2^{k_2}\sum_{j\geq
0}2^{-j/2}\|\mathbf{1}_{D_{k,\infty,j}}\cdot p(\xi,\mu)
\cdot (f_{k_1,j_1}\ast f_{k_2,j_2})\|_{L^2}\\
&\lesssim (k_2-k) 2^{k_2+k/2}\cdot 2^{j_1/2}\|p(\xi_1,\mu_1)\cdot
f_{k_1,j_1}\|_{L^2}\cdot 2^{j_2/2}\|p(\xi_2,\mu_2)\cdot
f_{k_2,j_2}\|_{L^2}.
\end{split}
\end{equation}

Instead of \eqref{pximu1} we now have
\begin{equation} \label{pximu2}
p(\xi,\mu) \lesssim 2^{2k_2 - k} (p(\xi_1,\mu_1)+p(\xi_2,\mu_2))
\end{equation}
which shows that the bound \eqref{on81} follows from
\eqref{jj1} for $j/2 \geq 2k_2 - k/2$. Therefore
it remains to prove that
\begin{equation}\label{on281a}
\begin{split}
&2^{-j/2}\|\mathbf{1}_{D_{k,\infty,j}}\cdot p(\xi,\mu)\cdot (f_{k_1,j_1}\ast f_{k_2,j_2})\|_{L^2}\\
&\lesssim 2^{-k/2}\cdot 2^{j_1/2}\|p(\xi_1,\mu_1)\cdot
f_{k_1,j_1}\|_{L^2}\cdot 2^{j_2/2}\|p(\xi_2,\mu_2)\cdot
f_{k_2,j_2}\|_{L^2}.
\end{split}
\end{equation}

We still have \eqref{etapexi}, but now this leads to
\[
p(\xi,\mu) \lesssim 2^{k_2}  p(\xi_1,\mu_1) + 2^{-k/2}
2^{\max(j_1,j_2,j)/2}
\]
Then \eqref{on281a} reduces to  \eqref{on81b} and \eqref{on81c},
which follow as before from \eqref{jj1}.
\end{proof}

Finally, we estimate low-frequency interactions.

\newtheorem{Lemmab6}[Lemmab2]{Lemma}
\begin{Lemmab6}\label{Lemmab6}
Assume $k,k_1,k_2\in(-\infty,100]\cap\Z$, $u_{k_1}\in F_{k_1}$, and $v_{k_2}\in F_{k_2}$.
Then
\begin{equation}\label{on380}
\|P_k(\partial_x(u_{k_1}v_{k_2}))\|_{N_k}\lesssim 2^{(k+k_1+k_2)/2} \cdot
\|u_{k_1}\|_{F_{k_1}}\cdot \|v_{k_2}\|_{F_{k_2}}.
\end{equation}
\end{Lemmab6}

\begin{proof}[Proof of Lemma \ref{Lemmab6}] As in the proof of Lemma
  \ref{Lemmab2}, using the definitions and the bounds \eqref{sp7.2}
  and \eqref{sb21}, it suffices to prove that if $j_1, j_2\in\Z_+$,
  and $f_{k_i,j_i}:\mathbb{R}^3\to\mathbb{R}_+$ are supported
  in $D_{k_i,\infty,j_i}$, $i=1,2$, then
\begin{equation}\label{on381}
\begin{split}
&2^k\sum_{j\geq 0}2^{-j/2}\|\mathbf{1}_{D_{k,\infty,j}}\cdot
p(\xi,\mu)
\cdot (f_{k_1,j_1}\ast f_{k_2,j_2})\|_{L^2}\\
&\lesssim 2^{(k+k_1+k_2)/2}\cdot 2^{j_1/2}\|p(\xi_1,\mu_1)\cdot
f_{k_1,j_1}\|_{L^2}\cdot 2^{j_2/2}\|p(\xi_2,\mu_2)\cdot
f_{k_2,j_2}\|_{L^2}.
\end{split}
\end{equation}
We may assume that $k_1\leq k_2$ (which forces $k\leq k_2+4$). We use the simple bound
\begin{equation*}
p(\xi,\mu) \lesssim 2^{k_2 -k}(p(\xi_1,\mu_1)+p(\xi_2,\mu_2)),
\end{equation*}
and \eqref{jj2lowk}. The bound \eqref{on381} follows.
\end{proof}

\section{Dyadic bilinear estimates, II}\label{dyadic3}

In this section we prove several dyadic bounds which are used in
the proof of Proposition \ref{Lemmab1} (b). We estimate first low-frequency interactions.

\newtheorem{Lemmah1}{Lemma}[section]
\begin{Lemmah1}\label{Lemmah1}
Assume $k,k_1,k_2\in\Z\cap(-\infty,100]$, $u_{k_1}\in
\overline{F}_{k_1}$, and $v_{k_2}\in F_{k_2}$. Then
\begin{equation}\label{hh1}
\|P_k(\partial_x(u_{k_1}v_{k_2}))\|_{\overline{N}_k}\lesssim
2^{3k/2+k_2/2} \cdot \|u_{k_1}\|_{\overline{F}_{k_1}}\cdot
\|v_{k_2}\|_{F_{k_2}}.
\end{equation}
\end{Lemmah1}

\begin{proof}[Proof of Lemma \ref{Lemmah1}] Using the definitions and the bounds \eqref{sp7.2} and \eqref{sb21}, it suffices to prove that if $j_1,
j_2\in\Z_+$, and $f_{k_i,j_i}:\mathbb{R}^3\to\mathbb{R}_+$
are supported in $D_{k_i,\infty j_i}$, $i=1,2$, then
\begin{equation}\label{hh2}
\begin{split}
&2^k\sum_{j\geq 0}2^{-j/2}\|\mathbf{1}_{D_{k,\infty,j}}\cdot (f_{k_1,j_1}\ast f_{k_2,j_2})\|_{L^2}\\
&\lesssim 2^{3k/2+k_2/2} \cdot
2^{j_1/2}\|f_{k_1,j_1}\|_{L^2}\cdot
2^{j_2/2}\|p(\xi_2,\mu_2)\cdot f_{k_2,j_2}\|_{L^2}.
\end{split}
\end{equation}
This is a direct consequence of \eqref{jj2lowk}.
\end{proof}

We estimate now $\mathrm{High}\times\mathrm{High}\to\mathrm{Low}$ interactions.

\newtheorem{Lemmah2}[Lemmah1]{Lemma}
\begin{Lemmah2}\label{Lemmah2}
Assume $k,k_1,k_2\in\Z$, $k_1,k_2\geq \max(k-10,20)$, $u_{k_1}\in
\overline{F}_{k_1}$, and $v_{k_2}\in F_{k_2}$. Then
\begin{equation}\label{hh20}
\|P_k(\partial_x(u_{k_1}v_{k_2}))\|_{\overline{N}_k}\lesssim
2^{(3k_2-3|k|)/4} \cdot \|u_{k_1}\|_{\overline{F}_{k_1}}\cdot
\|v_{k_2}\|_{F_{k_2}}.
\end{equation}
\end{Lemmah2}

\begin{proof}[Proof of Lemma \ref{Lemmah2}] As in the proof of Lemma \ref{Lemmab4}, using the definitions and the bounds
\eqref{sp7.2} and \eqref{sb21}, it suffices to prove that if $j_1,
j_2\geq k_2$, and $f_{k_i,j_i}:\mathbb{R}^3\to\mathbb{R}_+$
are supported in $D_{k_i,\infty,j_i}$, $i=1,2$, then
\begin{equation}\label{hh21}
\begin{split}
&2^k2^{k_2-k_+}\sum_{j\geq 0}2^{-j/2}\|\mathbf{1}_{D_{k,\infty,j}}\cdot (f_{k_1,j_1}\ast f_{k_2,j_2})\|_{L^2}\\
&\lesssim 2^{(3k_2-3|k|)/4} \cdot
2^{j_1/2}\|f_{k_1,j_1}\|_{L^2}\cdot
2^{j_2/2}\|p(\xi_2,\mu_2)\cdot f_{k_2,j_2}\|_{L^2}.
\end{split}
\end{equation}
Assume first that
\begin{equation*}
k\geq 0.
\end{equation*}
Then, using \eqref{jj1}, the left-hand side of \eqref{hh21} is
dominated by
\begin{equation*}
\begin{split}
&C2^{k_2}k_22^{-(2k_2+k)/2}2^{j_1/2}\|f_{k_1,j_1}\|_{L^2}\cdot 2^{j_2/2}\|f_{k_2,j_2}\|_{L^2},
\end{split}
\end{equation*}
which suffices for \eqref{hh21}. Assume now that
\begin{equation}\label{hh22}
k\leq 0\text{ and }k+2k_2\geq 0.
\end{equation}
Then, using \eqref{jj1}, the left-hand side of \eqref{hh21} is
dominated by
\begin{equation*}
\begin{split}
&C2^{k+k_2}k_22^{-(2k_2+k)/2}2^{j_1/2}\|f_{k_1,j_1}\|_{L^2}\cdot 2^{j_2/2}\|f_{k_2,j_2}\|_{L^2},
\end{split}
\end{equation*}
which suffices for \eqref{hh21} in view of \eqref{hh22}. Finally,
assume that
\begin{equation}\label{hh23}
k+2k_2\leq 0
\end{equation}
Then, using \eqref{jj1}, the left-hand side of \eqref{hh21} is
dominated by
\begin{equation*}
\begin{split}
C2^{k+k_2}2^{j_1/2}\|f_{k_1,j_1}\|_{L^2}\cdot
2^{j_2/2}\|f_{k_2,j_2}\|_{L^2},
\end{split}
\end{equation*}
which suffices for \eqref{hh21} in view of \eqref{hh23}.
\end{proof}

Finally, we estimate $\mathrm{Low}\times\mathrm{High}\to\mathrm{High}$ interactions.

\newtheorem{Lemmah3}[Lemmah1]{Lemma}
\begin{Lemmah3}\label{Lemmah3}
Assume $k,k_1,k_2\in\Z$, $k\geq 20$, $k_2\leq k-10$, $|k_1-k|\leq
4$, $u_{k_1}\in \overline{F}_{k_1}$, and $v_{k_2}\in F_{k_2}$.
Then
\begin{equation}\label{hh30}
\|P_k(\partial_x(u_{k_1}v_{k_2}))\|_{\overline{N}_k}\lesssim
2^{k_2} \cdot \|u_{k_1}\|_{\overline{F}_{k_1}}\cdot
\|v_{k_2}\|_{F_{k_2}}\text{ if }k_2 \leq 0,
\end{equation}
and
\begin{equation}\label{hh30a}
\|P_k(\partial_x(u_{k_1}v_{k_2}))\|_{\overline{N}_k}\lesssim
k_22^{-k_2/2} \cdot \|u_{k_1}\|_{\overline{F}_{k_1}}\cdot
\|v_{k_2}\|_{F_{k_2}}\text{ if }k_2\geq 1.
\end{equation}
\end{Lemmah3}

\begin{proof}[Proof of Lemma \ref{Lemmah3}] As in the proof of Lemma \ref{Lemmab2}, using the definitions and the bounds
\eqref{sp7.2} and \eqref{sb21}, it suffices to prove that if $j_1,
j_2\geq k$, and $f_{k_i,j_i}:\mathbb{R}^3\to\mathbb{R}_+$
are supported in $D_{k_i,\infty,j_i}$, $i=1,2$, then
\begin{equation}\label{hh31}
\begin{split}
&2^k\sum_{j\geq k}2^{-j/2}\|\mathbf{1}_{D_{k,\infty,j}}\cdot (f_{k_1,j_1}\ast f_{k_2,j_2})\|_{L^2}\\
&\lesssim c(k_2) \cdot
2^{j_1/2}\|f_{k_1,j_1}\|_{L^2}\cdot
2^{j_2/2}\|p(\xi_2,\mu_2)\cdot f_{k_2,j_2}\|_{L^2}.
\end{split}
\end{equation}
where $c(k_2) = 2^{k_2}$ for $k_2 \leq 0$, respectively
$c(k_2) =k_22^{-k_2/2}$  for $k_2 > 0$.

Consider first the case $k_2 \leq 0$.  Since $j,j_1,j_1 \geq k$, the
above bound is a direct consequence of \eqref{jj2lowk}.

If $k_2\geq 1$ then the high modulation case $j \geq k+4 k_2$
is obtained directly from \eqref{jj2highk}. Therefore it remains to prove that
\begin{equation}\label{hh31a}
\begin{split}
& \|\mathbf{1}_{D_{k,\infty,j}}\cdot (f_{k_1,j_1}
\ast f_{k_2,j_2})\|_{L^2}\\
& \lesssim 2^{-k-k_2/2} \cdot
2^{(j+j_1+j_2)/2}\|f_{k_1,j_1}\|_{L^2}\cdot
\|p(\xi_2,\mu_2)\cdot f_{k_2,j_2}\|_{L^2}.
\end{split}
\end{equation}
This follows from  \eqref{jj1}.

\end{proof}

\newtheorem{Lemmah4}[Lemmah1]{Lemma}
\begin{Lemmah4}\label{Lemmah4}
Assume $k,k_1,k_2\in\Z$, $k\geq 20$, $k_1\leq k-10$, $|k-k_2|\leq
4$, $u_{k_1}\in \overline{F}_{k_1}$, and $v_{k_2}\in F_{k_2}$.
Then
\begin{equation}\label{hh40}
\|P_k(\partial_x(u_{k_1}v_{k_2}))\|_{\overline{N}_k}\lesssim 2^{-k_1/2}k_2 \cdot \|u_{k_1}\|_{\overline{F}_{k_1}}\cdot
\|v_{k_2}\|_{F_{k_2}}.
\end{equation}
\end{Lemmah4}

\begin{proof}[Proof of Lemma \ref{Lemmah4}] As in the proof of Lemma \ref{Lemmab2}, using the definitions and the bounds
\eqref{sp7.2} and \eqref{sb21}, it suffices to prove that if $j_1,
j_2\geq 0$, and $f_{k_i,j_i}:\mathbb{R}^3\to\mathbb{R}_+$
are supported in $D_{k_i,\infty,j_i}$, $i=1,2$, then
\begin{equation}\label{hh41}
\begin{split}
&2^k\sum_{j\geq 0}2^{-j/2}\|\mathbf{1}_{D_{k,\infty,j}}\cdot (f_{k_1,j_1}\ast f_{k_2,j_2})\|_{L^2}\\
&\lesssim 2^{-k_1/2}k_2 \cdot
2^{j_1/2}\|f_{k_1,j_1}\|_{L^2}\cdot
2^{j_2/2}\|p(\xi_2,\mu_2)\cdot f_{k_2,j_2}\|_{L^2}.
\end{split}
\end{equation}
Using \eqref{jj1}, the left-hand side of \eqref{hh41} is dominated
by
\begin{equation*}
2^k k_2\cdot
2^{-(2k+k_1)/2}\|f_{k_1,j_1}\|_{L^2}\cdot\|f_{k_2,j_2}\|_{L^2},
\end{equation*}
which suffices.
\end{proof}

\end{document}